\documentclass[10pt]{amsart}
\usepackage{latexsym}
\usepackage{amsfonts}

\usepackage{amsmath,amsthm,amssymb}

\newtheorem{theor}{Theorem}

\newtheorem{thm}{Theorem}[section]
\newtheorem{lem}[thm]{Lemma}
\newtheorem{cor}[thm]{Corollary}

\newtheorem{prop}[thm]{Proposition}
\theoremstyle{definition}
\newtheorem{defn}[thm]{Definition}
\newtheorem{notation}[thm]{Notation}
\newtheorem{conv}[thm]{Convention}
\newtheorem{rem}[thm]{Remark}
\newtheorem{exmp}[thm]{Example}

\begin{document}

\author[Ilya Kapovich]{Ilya Kapovich}

\address{\tt Department of Mathematics, University of Illinois at
  Urbana-Champaign, 1409 West Green Street, Urbana, IL 61801, USA
  \newline http://www.math.uiuc.edu/\~{}kapovich/} \email{\tt
  kapovich@math.uiuc.edu}

\title[The frequency space of a free group]{The frequency space of a free group}
\thanks
{The author acknowledges support from the {\it Swiss National
Science Foundation}.}

\begin{abstract}
  We analyze the structure of the \emph{frequency space} $Q(F)$ of a
  nonabelian free group $F=F(a_1,\dots,a_k)$ consisting of all
  shift-invariant Borel probability measures on $\partial F$ and
  construct a natural action of $Out(F)$ on $Q(F)$. In particular we
  prove that for any outer automorphism $\phi$ of $F$ the
  \emph{conjugacy distortion spectrum} of $\phi$, consisting of all
  numbers $||\phi(w)||/||w||$, where $w$ is a nontrivial conjugacy
  class, is the intersection of $\mathbb Q$ and a closed subinterval of
  $\mathbb R$ with rational endpoints.
\end{abstract}


\maketitle


\section{Introduction}\label{intro}

Let $F=F(a_1,\dots, a_k)$ be a free group of rank $k>1$ with a fixed
free basis $X=\{a_1,\dots, a_k\}$. We identify the hyperbolic boundary
$\partial F$ with the set of all semi-infinite freely reduced words in
$F$. Let $T=T_X:\partial F\to \partial F$ be the \emph{standard shift
  operator}, which erases the first letter of every such semi-infinite
word over $X^{\pm 1}$.

The \emph{frequency space} $Q(F)=Q_{X}(F)$ can be defined as the space
of all $T$-invariant Borel probability measures on $\partial F$.
Objects of this sort appear naturally in ergodic theory, but in the
present paper we concentrate on the algebraic and combinatorial
structure of $Q(F)$ rather than on the more traditional analytical and
probabilistic questions.  Moreover, there is a natural correspondence
between $Q(F)$ and the space of ``geodesic currents'' on $F$, that is
the space of scalar equivalence classes of all $F$-invariant positive
Borel measures on $\partial^2 F:=\{(x,y)\in \partial F\times \partial
F: x\ne y\}$ (in many instances one considers ``unordered'' currents,
that is measures which are invariant both with respect to the
$F$-action and the flip map $(x,y)\mapsto (y,x)$). The space of
geodesic currents on a surface group, and, more generally on a
word-hyperbolic group, has been extensively studied by
Bonahon~\cite{Bo86,Bo88,Bo91} (whose work served as a substantial
source of inspiration for the present paper) and other authors and it
plays an important role in modern 3-manifold topology. We provide a
brief discussion of geodesic currents in Section~\ref{sect:currents} below.

In this paper we study a natural action of $Out(F)$ on $Q(F)$ by
homeomorphisms. The existence of
such an action was well understood in the surface group case and it
can in fact could be established via a correspondence between $Q(F)$
and geodesic currents, where it was first studied by Reiner Martin~\cite{Ma}. 
We provide a direct construction here which, as
we will see, also yields valuable new information about how this
action looks like in the natural ``frequency coordinates'' on $Q(F)$.

We first observe (this also follows
from Bonahon's work and was proved by Martin in the geodesic currents setting~\cite{Ma}) that there is a natural map $\alpha:C\to Q(F)$,
where $C$ is the set of nontrivial conjugacy classes in $F$, with \emph{dense}
image in $Q(F)$.  Also, $\alpha(C)=\alpha(C_0)$ where
$C_0$ is the set of root-free nontrivial conjugacy classes and
$\alpha$ is injective on $C_0$. It turns out that the natural action
of $Out(F)$ on $C_0$ translates into a continuous action on $Q(F)$.

We will briefly discuss the map $\alpha$ here, which will also
motivate the name ``frequency space'' for $Q(F)$. Let $w$ be a
nontrivial conjugacy class in $F$. We will think of $w$ as a
\emph{cyclic word} that is a cyclically reduced word written on a
circle in a clockwise direction without specifying a base-point.  We
denote the length of $w$ by $||w||$.  For any freely reduced word $v$
one can define in the obvious way the number $n_w(v)$ of occurrences
of $v$ in $w$. Namely, we count the number of positions on the circle,
starting from which it is possible to read the word $v$ going
clockwise along the circle (and wrapping around more than once, if
necessary). Thus $n_w(v)$ is defined for all $v\in F$, even for those
that are longer than $w$. We also define the \emph{frequency} $f_w(v)$
of $v$ in $w$ as $f_w(v)=n_w(v)/||w||$.  A cyclic word $w$ defines a
$T$-invariant probability measure $\mu_w=\alpha(w)$ on $\partial F$ as
follows.  If $v\in F$ is a nontrivial element, we set
\[
\mu_w(Cyl(v)):=f_w(v),
\]
where $Cyl(v)=Cyl_X(v)$ consists of all semi-infinite freely reduced words with
initial segment $v$.  It is not hard to see that $\mu_w$ is indeed a
$T$-invariant probability measure on $\partial F$.  Moreover, there is
another way to view the measure $\mu_w$.  Namely, for any cyclically
reduced word $v$ defining the cyclic word $w$ the point
$v^{\infty}=vvv\dots \in \partial F$ is $T$-periodic with the period
equal to the length of the deepest root $w_0$ of $w$. The $T$-orbit
$A(w)=\{T^nv^{\infty}| n\ge 0\}$ is independent of $v$ and consists of
$||w_0||$ points. It is not hard to see that $\mu_w$ is precisely the
discrete measure uniformly supported on $A(w)$:

\[
\mu_w=\frac{1}{\#A(w)} \sum_{x\in A(w)} \delta_x.
\]

It follows that that $\alpha(C)$ is dense in $Q(F)$, since the finite
sets $A(w)$ are exactly the $T$-periodic orbits in $\partial F$. Thus
for any measure $\mu\in Q(F)$ one may call the number
$f_{\mu}(v):=\mu(Cyl(v))$ the \emph{frequency} of $v$ in $\mu$. Note
that the collection of frequencies $(f_{\mu}(v))_{v\in F}$ uniquely
determines $\mu\in Q(F)$ and thus $(f_{.}(v))_{v\in F}$ can be
considered as a global ``coordinate system'' for $Q(F)$.  
We prove (see Theorem~\ref{action} below) that the canonical action of
$Out(F)$ on $C_0$ extends  via $\alpha$ (uniquely since $\alpha(C_0)$
is dense in $Q(F)$) to an action of $Out(F)$ on $Q(F)$ by homeomorphisms.
Moreover, it turns out that in the frequency coordinates described
above the action of each $\phi\in Out(F)$ on $Q(F)$ is expressible by
\emph{fractional linear transformations}, where each coordinate
function $f_{\phi(.)}(v)$ (with $v\in F$ being a fixed element)
depends in a fractional-linear manner on only finitely many frequency
coordinates $(f_{.}(v))_{v\in F}$ of the argument. Moreover, these
transformations are, in a sense, linear. Namely, for a fixed point
$p\in Q(F)$ all the fractional-linear coordinate functions
$f_{\phi(p)}(v)$, $v\in F$, have the same denominators.

We also study in detail the structure of $Q(F)$ and the image of the
map $\alpha$, concentrating on the description of $Q(F)$ as the
inverse limit of a sequence of affine maps between finite-dimension
convex compact polyhedra $Q_m$. Roughly speaking $Q_m$ captures
``finitary chunks'' of $T$-invariant measures $\mu$, namely
$(f_{\mu}(u))_{|u|=m}$. In particular, we obtain an explicit geometric
and algorithmically verifiable criterion (Theorem~\ref{realize})
characterizing those tuples of rational numbers $q=(q_u)_{u\in F,
  |u|=m}$ that can be realized as frequency tuples of cyclic words.

Although the definition of $Q(F)=Q_X(F)$ is less equivariant
than that of the space of geodesic currents (because of the dependence
on $X$), it turns out that the action of $Out(F)$ on $Q(F)$ is
particularly informative for algebraic purposes, as it contains some
new interesting algebraic, geometric and algorithmic information about
automorphisms of $F$.  In fact, it appears that the action of $Out(F)$
on $Q(F)$ often provides data of a different kind from what one
usually gets by traditional outer space and train-track methods.

Our main applications concern the metric distortion properties of
automorphisms. For an automorphism $\phi\in Aut(F)$ we define the
\emph{conjugacy distortion spectrum} $I(\phi)$ of $\phi$  (with
respect to a fixed basis $A$ of $F$) as
\[
I(\phi):=\{\frac{||\phi(w)||}{||w||}: w\in F, w\ne 1\}.
\]
It is obvious that $I(\phi)\subseteq \mathbb Q$ and one can easily
show that there exist numbers $0<C_1,C_2<\infty$ such that
$I(\phi)\subseteq [C_1,C_2]$. Surprisingly, it is possible to obtain a
very precise description of the possible shape of $I(\phi)$. Thus we prove:

\begin{theor}\label{A}
For any $\phi\in Aut(F)$ there exists a closed interval $J\subseteq \mathbb
R$ such that

\[
I(\phi)=J\cap \mathbb Q.
\]
Moreover, the endpoints of $J$ are positive rational numbers that are
algorithmically computable in terms of $\phi$.
\end{theor}

Theorem~\ref{A} implies, in particular, that the ``extremal distortion
factors'' $\sup I(\phi)$ and $\inf I(\phi)$ actually belong to $I(\phi)$ and are
therefore realized as distortion factors of some nontrivial conjugacy classes in $F$.

In a subsequent work with Kaimanovich and Schupp~\cite{KSS} we combine
the results of this paper with a number of algebraic and probabilistic
methods to obtain certain rigidity results regarding $I(\phi)$. In
particular, it turns out that for any $\phi\in Out(F)$ there exists a
nontrivial conjugacy class $w$ such that $||w||=||\phi(w)||$.
Moreover, there exists an open interval $I_0\subseteq \mathbb R$
containing $1$ such that for any $\phi\in Out(F)$ either $\phi$ is
equal in $Out(F)$ to an automorphism induced by a permutation of
$X^{\pm 1}$, in which case $I(\phi)=\{1\}$, or $\mathbb Q\cap
I_0\subseteq I(\phi)$.

We also obtain some applications  to \emph{hyperbolic}
automorphisms of free groups. An outer automorphism $\phi\in Out(F)$
is \emph{strictly hyperbolic} if there is $\lambda>1$ such that
\[
\lambda ||w||\le \max\{ ||\phi(w)||, ||\phi^{-1}||\} \text{ for every
  cyclic word } w.
\]
An outer automorphism $\phi\in Out(F)$ is \emph{hyperbolic} if there
is $n>0$ such that $\phi^n$ is strictly hyperbolic (this definition is
equivalent to the one used in \cite{BF92}, as observed in \cite{BFH}).

This notion is important because of its connection with the
Combination Theorem of Bestvina and Feign~\cite{BF92,BF96} (see also
\cite{Ga}) and the structure of free-by-cyclic groups. Peter
Brinkmann~\cite{Br} proved that $\phi\in Out(F)$ is hyperbolic if and
only if $\phi$ does not have nontrivial periodic conjugacy classes.
Using the Combination Theorem this implies that for $\psi\in Aut(F)$
the free-by-cyclic group $G=F\rtimes_{\psi} {\mathbb Z}$ is hyperbolic
if and only if the outer automorphism $\phi$ defined by $\psi$ has no
periodic conjugacy classes, that is if and only if $G$ does not
contain ${\mathbb Z}\times {\mathbb Z}$-subgroups.

For $\phi\in Out(F)$ let $\lambda_0(\phi)$ be the infimum over all
nontrivial cyclic words $w$ of \[\max\{ \frac{||\phi(w)||}{||w||},
\frac{||\phi^{-1}(w)||}{||w||}\}.\]

\begin{theor}\label{B}
  Let $\phi\in Out(F)$. Then $\lambda_0(\phi)$ is a rational number,
  algorithmically computable in terms of $\phi$. The outer
  automorphism $\phi$ is strictly hyperbolic if and only if
  $\lambda_0(\phi)>1$, and hence it is algorithmically decidable
  whether $\phi$ is strictly hyperbolic.
\end{theor}

This provides an algorithm (apparently the first known one) for
deciding if $\phi$ is strictly hyperbolic: compute $\lambda_0(\phi)$
and check if $\lambda_0(\phi)>1$. Moreover, a detailed analysis of the
proof produces an explicit worst-case complexity bound for the above
algorithms, namely, double exponential time in terms of the length of
$\phi$ as a product of generators of $Out(F)$. As we note later, one
can also use Theorem~\ref{B}, together with a result of Brinkmann~\cite{Br}, to provide a new algorithm for deciding
if an automorphism is hyperbolic.

The author is very grateful to Peter Brinkmann, Paul Schupp, Douglas
West, Igor Rivin, Ilia Binder, Vadim Kaimanovich, Tatiana
Smirnova-Nagnibeda, Martin Lustig and Kim Whittlesey
for numerous helpful comments, suggestions and discussions. In
particular the current proof of Lemma~\ref{real1} is due to Brinkmann
and it is much shorter than the author's original proof. The idea of
the proof of
Theorem~\ref{realize} was suggested to the author by Igor Rivin. The
author also thanks Pierre de la Harpe and the Swiss National Science
Foundation for providing financial support for a visit by the author
to the University of Gen\'{e}ve in June 2003, where much of the work
on this paper was conducted.

After the first version of this paper has been released, the author
was alerted by Mladen Bestvina of the existence of an earlier
unpublished dissertation of his UCLA student Reiner Martin~\cite{Ma}.
In his 1995 PhD thesis~\cite{Ma} Martin studied, in particular, the
space of geodesic currents on a free group. In particular, he obtained
a description of an action of $Out(F)$ on the space of currents which,
via its identification with the frequency space, can be seen to
coincide with the action constructed here. Moreover, there are,
inevitably, many similarities in the general treatment of the space of
currents by Martin and frequency space in the present paper. However,
as it turned out, our approaches also differ in many details and in
most applications.  Thus our treatment of $Q_X(F)$ as being
approximated by finite-dimensional polyhedra $Q_m$ is new, as are our
main applications, namely the ``realization theorem''
Theorem~\ref{realize}, Theorem~\ref{A} and Theorem~\ref{B}.  On
the other hand, Martin's thesis contains many interesting results that
are not addressed in the present paper.

The author is especially grateful to Mladen Bestvina for letting him
know about Martin's work and providing a copy of Martin's dissertation. 

\section{Frequencies in cyclic words}

\begin{conv}
  
  Fix an integer $k\ge 2$ and $F=F(a_1,\dots, a_k)$. Denote
  $X=\{a_1,\dots, a_k\}^{\pm 1}$. We denote by $CR$ the set of all
  cyclically reduced words in $F$.
  
  A \emph{cyclic word} is an equivalence class of cyclically reduced
  words, where two cyclically reduced words are equivalent if they are
  cyclic permutations of each other. If $u$ is a cyclically reduced
  word, we denote by $(u)$ the cyclic word defined by $u$. If $u$ is a
  freely reduced word, we denote the length of $u$ by $|u|$ and the
  length of the cyclically reduced form of $u$ by $||u||$. If $w=(v)$
  is a cyclic word, we denote by $||w||:=||v||$ the \emph{length} of
  $w$.

  We will denote the set of all nontrivial cyclic words by $C$.
  
  Let $w$ be a nontrivial cyclic word and let $u$ be a nontrivial
  freely reduced word. We denote by $n_w(u)$ the number of occurrences
  of $u$ in $w$.  This notation has obvious meaning even if $|u|>|w|$.
  Namely, let $w=(z)$ for a cyclically reduced word $z$. Take the
  smallest $p>0$ such that $|z^{p-1}|\ge 2|u|$ and count the number of
  those $i, 0\le i<|w|$ such that $z^p\equiv z_1uz_2$ where $|z_1|=i$.
  This number by definition is $n_w(u)$. If $u=1$ we define
  $n_w(u):=||w||$.
  
  Another way to think of $n_w(u)$ is as follows. Identify a cyclic
  word $w$ with a circle subdivided in $||w||$ directed edges labeled
  by letters from $X$ so that starting at any position and going
  clockwise around the circle once we read a cyclically reduced word
  defining $w$. Then $n_w(u)$ is the number of vertices on the circle
  starting from which and going clockwise it is possible to read the
  word $u$ without getting off the circle.

  Also, we denote $f_w(u):=\frac{n_w(u)}{||w||}$ and call it the
  \emph{frequency} of $u$ in $w$.
  
  Thus, for example, $n_a(aaa)=1$ and $f_w(1)=1$ for any $w$.

  If $w$ is a nontrivial freely reduced word (that is not necessarily
  cyclically reduced), and $v$ is a freely reduced word, we denote by
  $n_w(v)$ the number of occurrences of $v$ in $w$. If $|w|=n>0$ then
  by definition $n_w(v)$ is the number of those $i, 0\le i<n$ for
  which $w$ decomposes as a freely reduced product $w=w'vw''$ with
  $|w'|=i$. Thus if $|v|\le |w|$ then necessarily $n_w(v)=0$ (unlike
  the situation if $w$ is a cyclic word).
\end{conv}

\begin{lem}\label{red}\cite{KKS}
  Let $w$ be a nontrivial cyclic word. Then:

\begin{enumerate}
  
\item For any $m\ge 0$ and for any freely reduced word $u$ with
  $|u|=m$ we have:
\[
n_w(u)=\sum \{ n_w(ux): x\in X, |ux|=|u|+1\}=\sum \{ n_w(xu): x\in X,
|xu|=|u|+1\},
\]
and
\[
f_w(u)=\sum \{ f_w(ux): x\in X, |ux|=|u|+1\} =\sum \{ f_w(xu): x\in X,
|xu|=|u|+1\}.
\]

\item For any $m\ge 1$
\[
\sum_{|u|=m} n_w(u)=||w|| \quad \text{ and } \sum_{|u|=m} f_w(u)=1.
\]

\item For any $s>0$ and any $u\in F$
\[
n_{w^s}(u)=sn_w(u) \quad \text{ and }\quad f_{w^s}(u)=f_w(u).
\]
\end{enumerate}

\end{lem}

\begin{proof}
  Parts (1) and (3) are obvious. We establish (2) by induction on $m$.
  For $m=1$ the statement is clear. Suppose that $m>1$ and that (2)
  has been established for $m-1$.
  
  We have:

\[
\sum_{|u|=m} n_w(u) = \sum_{\overset{|v|=m-1, x\in X:} {|vx|=m}}
n_w(vx) = \sum_{|v|=m-1} n_w(v)=||w||,
\]
as required.

\end{proof}

Note that $0\le f_w(v)\le 1$ and $f_w(v)\in \mathbb Q$.  Part (2) of
Lemma~\ref{red} implies that if $w$ is a cyclic word and $m\ge 1$ is
an integer, then the set of frequencies $(f_w(v))_{\{v\in F: |v|=m\}}$
defines a point in the standard simplex of dimension $2k
(2k-1)^{m-1}-1$ standardly embedded in the Euclidean space of
dimension $2k(2k-1)^{m-1}$.

\begin{conv}
  Denote $D(m):=2k(2k-1)^{m-1}$ for $m>0$ and $D(0):=0$. Note that for
  any $m\ge 0$ $D(m)$ is exactly the number of elements of length $m$
  in $F$. We also identify ${\mathbb R}^{0}={\mathbb R}^{D(0)}=\{1\}$.
  
  In the subsequent discussion we will think of the coordinates of
  points $q\in {\mathbb R}^{D(m)}$ as indexed by elements $v\in F$
  with $|v|=m$.  Thus a point $q\in {\mathbb R}^{D(m)}$ is a
  $D(m)$-tuple $(q_v)_{|v|=m}$ of real numbers.
\end{conv}

\begin{defn}\label{defn:qm}
  Let $m\ge 0$ be an integer.
  We define the set $Q_m\subseteq {\mathbb R}^{D(m)}$ as the
  collection of all points $q=(q_v)_{|v|=m}\in {\mathbb R}^{D(m)}$
  that satisfy the following conditions:

\begin{enumerate}
\item We have
$\displaystyle\sum_{|v|=m} q_v =1$.

\item For each $v\in F$ with $|v|=m$ we have $q_v\ge 0$.

\item For each $u\in F$ with $|u|=m-1$ we have

\[
\sum_{\{x\in X:|ux|=m\}} q_{ux} =\sum_{\{y\in X:|yu|=m\}} q_{yu}.
\]
We also define $Q_0:=\{1\}$. Sometimes, we will also denote $Q_m$ as
$Q_{m,k}$ to stress the dependence of this notion on the rank $k$ of
the free group $F$.
\end{enumerate}

\end{defn}
Thus $Q_m$ is a convex compact polyhedron in ${\mathbb R}^{D(m)}$
that, for $m>0$, is contained in the standardly embedded
$D(m)-1$-dimensional simplex in ${\mathbb R}^{D(m)}$.

Note that for $m=1$ condition (3) in the above definition is vacuously
satisfied.

\begin{defn}\label{defn'}
  Let $m\ge 2$. We define maps $\hat \pi_m: {\mathbb R}^{D(m)}\to
  {\mathbb R}^{D(m-1)}$ as follows. Let $q=(q_v)_{|v|=m}\in {\mathbb
    R}^{D(m)}$. Then for any $u\in F$ with $|u|=m-1$ we set the $u$s
  coordinate of $\hat \pi_m(q)$ to be:

\[
[\hat \pi_m(q)]_u:=\sum_{\{x\in X:|ux|=m\}} q_{ux}.
\]
\end{defn}

\begin{lem}\label{lem:image}
  For any $m\ge 1$ we have
\[\hat \pi_m(Q_m)\subseteq Q_{m-1}\]
\end{lem}

\begin{proof}
  
  Let $q\in Q_m$ and $z=\hat \pi_m(q)$. It is obvious that for every
  $u\in F$ with $|u|=m-1$ $z_u\ge 0$.
  
  We have
\[
\sum_{|u|=m-1} z_u = \sum_{|u|=m-1} \sum_{\{x\in X:|ux|=m\}} q_{ux}=
\sum_{|v|=m} q_v =1,
\]
so that condition (1) of Definition~\ref{defn:qm} holds for $z$.

To verify condition (3) of Definition~\ref{defn:qm}, let $u'\in F$ be
an element with $|u'|=m-2$.

Then

\begin{align*}
  &\sum_{\{x\in X:|u'x|=m-1\}} z_{u'x} = \sum_{\{x\in X:|u'x|=m-1\}}
  \sum_{\{x_1\in X:|u'xx_1|=m\}} q_{u'xx_1}=\\
  &\sum_{\{x\in X:|u'x|=m-1\}} \sum_{\{y\in X:|yu'x|=m\}}
  q_{yu'x}=\sum_{\{x,y\in X: |yu'x|=m\}}q_{yu'x}.
\end{align*}

Here the first equality holds by the definition of $\hat \pi_m$ and
the second equality holds by condition (3) of Definition~\ref{defn:qm}
for $q$.

Similarly,

\begin{align*}
  &\sum_{\{y\in X:|yu'|=m-1\}} z_{yu'} = \sum_{\{y\in X:|yu'|=m-1\}}
  \sum_{\{x\in X:|yu'x|=m\}} q_{yu'x}=\\
  &=\sum_{\{x,y\in X: |yu'x|=m\}}q_{yu'x}.
\end{align*}

Thus
\[
\sum_{\{x\in X:|u'x|=m-1\}} z_{u'x}=\sum_{\{y\in X:|yu'|=m-1\}}
z_{yu'}
\]
and condition (3) of Definition~\ref{defn:qm} holds for $z$.
\end{proof}

\begin{defn}\label{defn:pm}
  In view of Lemma~\ref{lem:image} for $m\ge 2$ we define
\[
\pi_m: Q_m\to Q_{m-1}
\]
to be the restriction of $\hat \pi_m$ to $Q_m$. For $m=1$ define
$\pi_1:Q_1\to Q_0$ as $\pi_1(q)=1$ for each $q\in Q_1$. Note that this
agrees with the formula defining $\hat q_m$ in Definition~\ref{defn'},
since for any $q\in Q_1$ $\sum_{|x|=1} q_x=1$.
\end{defn}

We will later on see that the map $\pi_m$ is actually ``onto''.

Any cyclic word $w$ comes equipped with its set of frequencies
$(f_w(v))_{|v|=m}$ which by Lemma~\ref{red} satisfy all the conditions
of Definition~\ref{defn:qm}.

\begin{defn}
  Let $m\ge 1$. We define the map $\alpha_m: C\to Q_m$ as follows.
  For a cyclic word $w\in C$ and $v\in F$ with $|v|=m$ we set
\[
[\alpha_m(w)]_v:= f_w(v).
\]
\end{defn}

We have seen that the sets $Q_m$ are compact convex Euclidean
polyhedra and that the maps $\pi_m: Q_m\to Q_{m-1}$ are affine for
each $m\ge 1$. We can now formulate the main definition of the paper.

\begin{defn}[The Frequency Space]
  Consider the systems of maps
\[
\dots \to Q_m \overset{\pi_m}{\to} Q_{m-1}\overset{\pi_{m-1}}{\to}
\dots \to Q_2\overset{\pi_2}{\to} Q_1 \overset{\pi_1}{\to}
Q_0\tag{$\dag$}
\]

We define the \emph{frequency space} $Q(F)=Q_X(F)$ to be the inverse
limit of this sequence:
\[
Q(F):=\underset{\longleftarrow}{\lim} (Q_m, \pi_m)
\]
\end{defn}

Thus an element of $Q(F)$ is a tuple of numbers $(q_v)_{v\in F}$ that
are consistent with respect to the sequence $(\dag)$.  This means
that:

(a) For any $m\ge 1$ $(q_v)_{|v|=m}$ defines a point $q_m$ of $Q_m$.

(b) For any $m\ge 2$ $\pi_m(q_m)=q_{m-1}$.

Note that $Q(F)$ comes equipped with a canonical inverse limit
topology induced by the inclusion $Q(F) \subseteq \prod_{m\ge 0} Q_m$.
Namely, two points $q=(q_v)_{v\in F}$ and $q'=(q'_v)_{v\in F}$ are
close if there is large $M\ge 1$ such that $q_m$ is close to $q_m'$ in
$Q_m$ for all $1\le m\le M$ (or, equivalently, that $q_M$ is close to
$q_M'$ in $Q_M$).

\begin{rem}
  It is easy to see that the natural maps $\alpha_m:C\to Q_{m}$ define
  a map
\[
\alpha: C\to Q(F)
\]

If $w\in C$ is a cyclic word, then its image $\alpha(w)\in Q(F)$ can
be thought of as an infinite tuple $(f_w(v))_{v\in F}$.
\end{rem}

Part~(3) of Lemma~\ref{red} immediately implies that the maps
$\alpha_m$ and $\alpha$ have a certain "projective" character in the
sense that the image of $w$ depends only on the semigroup generated by
$w$:

\begin{lem}\label{power}
  Let $w\in C$ be a nontrivial cyclic word. Then for any $s\ge 1$ we
  have $\alpha(w)=\alpha(w^s)$ and $\alpha_m(w)=\alpha_m(w^s)$, where
  $m\ge 0$.
\end{lem}

Moreover, having common positive powers is the only reason why two
cyclic words can have the same image in $Q(F)$ under $\alpha$:

\begin{prop}\label{imalpha}
  Let $w, u$ be nontrivial cyclic words. Then $\alpha(w)=\alpha(u)$ if
  and only if both $w$ and $u$ are positive powers of the same cyclic
  word.
\end{prop}

\begin{proof}
  The "if" implication is obvious by Lemma~\ref{power}. Suppose now
  that $\alpha(w)=\alpha(u)$. Hence $\alpha_m(w)=\alpha_m(u)$ for
  every $m\ge 0$.

  There exist some $s,t>0$ such that $||w^s||=||u^t||$. Denote
  $j=||w^s||=||u^t||$.  Then $\alpha_m(w^s)=\alpha_m(u^t)$ for each
  $m\ge 0$ and, in particular $\alpha_j(w^s)=\alpha_j(u^t)$.
  
  Since $j=||w^s||$, whenever $f_{w^s}(v)>0$ for $v\in F$ with
  $|v|=j$, then $v$ must be a cyclically reduced word with $(v)=w^s$.
  The same is true for $u^t$. Choose $v\in F$ with $|v|=j$ such that
  $f_{w^s}(v)=f_{u^t}(v)>0$. Then $w^s=(v)=u^t$, which implies the
  statement of Proposition~\ref{imalpha}.
\end{proof}

\section{The frequency space as the space of invariant measures}

Recall that if $(\Omega, {\mathcal F}, \mu)$ is a probability space,
and $T:\Omega\to \Omega$ be a measurable map, then $\mu$ is said to be
\emph{$T$-invariant} if for any $A\in \mathcal F$ we have
$\mu(A)=\mu(T^{-1}A)$.

\begin{conv}
  Suppose $F=F(a_1,\dots, a_k)$ if a free group of finite rank $k>1$.
  Then the hyperbolic boundary $\partial F$ is naturally identified
  with the set of semi-infinite freely reduced words in $\{a_1,\dots,
  a_k\}^{\pm 1}$ corresponding to geodesic rays from $1$ in the
  standard Cayley graph of $F$. We endow $\partial F$ with the
  $\sigma$-algebra $\mathcal F$ of Borel sets. If $u$ is a freely
  reduced word of length $n$ then by $Cyl(u)$ we denote the set of all
  $x\in \partial X$ that have $u$ as initial segment.  Denote by
  $T:\partial F\to \partial F$ the \emph{shift operator} which erases
  the first letter of every semi-infinite freely reduced word.
\end{conv}

\begin{thm}\label{thm:ident}
  The frequency space $Q(F)$ is canonically identified with the space
  of all $T$-invariant Borel probability measures on $\partial F$.
\end{thm}

\begin{proof}
  Indeed, suppose $\mu$ is a $T$-invariant Borel probability measure
  on $\partial F$.
  
  For each $v\in F$ put $q_v:=\mu(Cyl(v))$. We claim that the tuple
  $(q_v)_{v\in F}$ defines a point of $Q(F)$. We need to show that for
  any $m\ge 1$ $q_m:=(q_v)_{|v|=m}\in Q_m$ and that for $m\ge 1$ we
  have $\pi_m(q_m)=q_{m-1}$.
  
  For any $m\ge 0$

\[
1=\mu(\partial F)=\sum_{|v|=m} \mu(Cyl(v))=\sum_{|v|=m} q_v,
\]

so that condition (1) of Definition~\ref{defn:qm} holds for $q_m$.
For any $v\in F$ we have
\[
Cyl(v)=\bigsqcup_{\{x\in X: |vx|=|v|+1\}} Cyl(vx)
\]
which implies
\[
q_v=\sum_{\{x\in X: |vx|=|v|+1\}} q_{vx}.
\]

Moreover, $T^{-1}(Cyl(v))=\bigsqcup_{\{x\in X: |xv|=|v|+1\}} Cyl(xv)$
and therefore by $T$-invariance of $\mu$
\[
q_v=\sum_{\{x\in X: |xv|=|v|+1\}} q_{xv}
\]
This implies that for each $m\ge 1$ all the conditions of
Definition~\ref{defn:qm} hold for $q_m$ and so $q_m\in Q_m$.

Moreover, by definition of $\pi_m$ the above equation also implies
that for any $m\ge 1$ we have $\pi_m(q_m)=q_{m-1}$, as required.  Thus
indeed $q=(q_v)_{v\in F}\in Q(F)$.

Suppose now that $q=(q_v)_{v\in F}\in Q(F)$.  Thus for any $m\ge 1$
$q_m:=(q_v)_{|v|=m}\in Q_m$ and $\pi_m(q_m)=q_{m-1}$ provided $m\ge
1$.  We will specify a Borel measure $\mu$ on $\partial F$ by its
values on all the sets $Cyl(v), v\in F$.  For $v\in F$ set

\[
\mu(Cyl(v)):=q_v.
\]

Then, by conditions (1)-(3) of Definition~\ref{defn:qm} and because
$\pi_m(q_m)=q_{m-1}$ it is easy to see that $\mu$ is indeed a
$T$-invariant Borel probability measure on $\partial F$.

We have constructed two maps between $Q(F)$ and the set of
$T$-invariant Borel probability measures on $\partial F$ that are
mutually inverse. This maps provide the required identification
asserted in Theorem~\ref{thm:ident}. Moreover, it is easy to see that the
inverse limit topology on $Q(F)$ coincides with the weak topology on
the space of Borel $T$-invariant measures.
\end{proof}

Recall that for a cyclic word $w$ the measure $\alpha(w)$ on $\partial
F$ was defined by $\alpha(w)(Cyl(u))=f_w(u)$ for $u\in F$. There is
another simple description of $\alpha(w)$.

If $v$ is a cyclically reduced word, we denote by $v^{\infty}$ the
element of $\partial F$ corresponding to the geodesic ray from $1$
with the label $vvvv\dots$. For a cyclic word $w$ denote
\[
A(w):=\{ v^{\infty}| (v)=w\}.
\]
It is not hard to see that $A(w)$ is a finite subset of $\partial F$
consisting of $||w_0||$ elements where $w=w_0^s$ with maximal possible
$s\ge 1$. In particular, if $w$ is not a proper power then
$\#A(w)=||w||$.  Recall that a point $x\in \partial F$ is said to be
$T$-\emph{periodic} if there is $n>0$ such that $T^nx=x$.

The following is an easy corollary of the definitions:

\begin{lem}\label{lem:a}
  The following hold:
\begin{enumerate}
\item Let $w$ be a cyclic word that is not a proper power and let
  $u\in F$ be a nontrivial freely reduced word. Then $w_u$ is equal to
  the number of those $x\in A(w)$ that have initial segment $u$, that
  is, to the number of those $x\in A(w)$ that are contained in $C(u)$.
\item A point $x\in \partial F$ is $T$-periodic if and only if
  $x=v^{\infty}$ for some cyclically reduced word $v$.
\item For any cyclic word $w$ and any $v\in F$ with $(v)=w$ the
  $T$-orbit of $v^{\infty}$ is precisely $A(w)$:
\[
A(w)=\{T^n v^{\infty} | n\ge 0\}.
\]
\end{enumerate}
\end{lem}

This in turn immediately implies:

\begin{lem}\label{lem:b}
  For any cyclic word $w$ the measure $\alpha(w)$ is the discrete
  measure uniformly supported on $A(w)$:
\[
\alpha(w)=\frac{1}{\#A(w)} \sum_{x\in A(w)} \delta_x.
\]
\end{lem}

\begin{cor}\label{cor:dense}
  The set $\alpha(C)$ is dense in $Q(F)$.
\end{cor}
\begin{proof}
  By Theorem~\ref{thm:ident} the space $Q(F)$ is exactly the set of
  $T$-invariant Borel probability measures on $\partial F$.  It
  follows from the basic results of symbolic dynamics~\cite{KH}
  (cf~\cite{Par,Ox}) that the set $Z$ of discrete measures uniformly
  supported on $T$-periodic orbits is dense in $Q(F)$. By
  Lemma~\ref{lem:a} and Lemma~\ref{lem:b} $Z=\alpha(C)$ and the
  statement follows.
\end{proof}

\section{Analyzing the image of the map $\alpha_m$}

Our next goal is to understand the image of the map $\alpha_m$, namely
$\alpha_m(C)\subseteq Q_m$. Obviously, each point in $\alpha_m(C)$
must have all rational coordinates, but it is not clear for the moment
if this conditions is also sufficient.

To study this question we need to introduce the following useful
notion.

\begin{defn}[Initial graph]\label{defn:G}
  Let $m\ge 2$ and let $q=(q_v)_{|v|=m}$ be a point in $Q_m$. We
  define the \emph{initial graph} $\Gamma_q$ as follows.
  
  The vertex set of $\Gamma_q$ is the set $\{u\in F : |u|=m-1\}$.
  
  For each $v\in F$ with $|v|=m-1$ there is a directed edge in
  $\Gamma_q$ with label $q_v$ from the vertex $v_{-}$ to the vertex
  $v_{+}$.
  
  Here $v_{-}$ denotes the initial segment of $v$ of length $m-1$ and
  $v_{+}$ denotes the terminal segment of $v$ of length $m-1$.  Thus
  $\Gamma_q$ is a labelled directed graph with the sum of the
  edge-labels equal to $1$.

  We also denote by $\Gamma_q'$ the graph obtained from $\Gamma_q$ by
  first removing all the edges labeled by $0$ and then removing all
  isolated vertices from the result. We call $\Gamma_q'$ the
  \emph{improved initial graph} of $q$.  The edge-set of $\Gamma_q'$
  is in 1-to-1 correspondence with the set of those $v\in F, |v|=m$
  for which $q_v>0$.
  
  If $w$ is a nontrivial cyclic word and $q=\alpha_m(w)$, we denote
  the graph $\Gamma_q$ by $\Gamma_{w,m}$ (or just $\Gamma_w$ if the
  value of $m$ is fixed). Similarly, we denote $\Gamma_q'$ by
  $\Gamma_{w,m}'$ (or just $\Gamma_w'$).
\end{defn}

The following lemma is an immediate corollary of
Definition~\ref{defn:qm}.

\begin{lem}
  If $m\ge 2$ then for any $q\in Q_m$ the graphs $\Gamma_q$ and
  $\Gamma_q'$ have the following properties:

\begin{enumerate}
\item The sum of the edge-labels in each of $\Gamma_q, \Gamma_q'$ is
  equal to $1$.
  
\item For any vertex $u$ of $\Gamma_q$ the in-degree of $u$ is equal
  to the out-degree of $u$ in $\Gamma_q$. The same is true for
  $\Gamma_q'$. Moreover, for any vertex $u$ of $\Gamma_q$ its
  in-degree in $\Gamma_q'$ is the same as its in-degree in
  $\Gamma_q'$.
\end{enumerate}

\end{lem}

\begin{notation}
  Let $\Gamma$ be a labeled directed graph and let $r\in \mathbb R$.
  We denote by $r\Gamma$ the graph obtained from $\Gamma$ by
  multiplying the label of every edge of $\Gamma$ by $r$.
  
  Let $\Gamma$ be a labeled directed graph where the label of any edge
  is a positive integer. We denote by $[\Gamma]$ the directed
  unlabelled graph obtained from $\Gamma$ by replacing every directed
  edge $e$ from a vertex $u$ to a vertex $v$ with label $M>0$ by $M$
  directed unlabelled edges from $u$ to $v$.

  If $\Gamma$ is a directed unlabelled graph, and $n\ge 1$ is an
  integer, we denote by $n\Gamma$ the directed unlabelled graph
  obtained from $\Gamma$ by multiplying the number of directed edges
  between any two vertices of $\Gamma$ by $n$.
\end{notation}

By a \emph{labelled directed graph} we will mean a digraph where every
edge is equipped with a real number called the \emph{label} of the
edge. A \emph{cycle} in a directed graph is a closed directed
edge-path. A \emph{circuit} in a directed graph is an equivalence
class of cycles, where two cycles are equivalent if they are cyclic
permutations of each other. Thus the notion of a circuit is similar to
that of a cyclic word in the free group context. An \emph{Euler cycle}
in a directed graph is a cycle that passes through each edge of the
graph exactly once. An \emph{Euler circuit} in a directed graph is an
equivalence class of an Euler cycle.

We need the following elementary fact from graph theory.

\begin{lem} Let $\Gamma$ be a directed graph.
\begin{enumerate}
\item Suppose the underlying undirected graph of $\Gamma$ is
  connected. Then $\Gamma$ possesses an Euler circuit (in the digraph
  sense) if and only if for every vertex $u$ of $\Gamma$ the in-degree
  at $u$ is equal to the out-degree at $u$.
  
\item Suppose that for every vertex $u$ of $\Gamma$ the in-degree at
  $u$ is equal to the out-degree at $u$. Then the underlying
  undirected graph of $\Gamma$ is connected if and only if $\Gamma$ is
  connected as a digraph (that is for any vertices $u_1$, $u_2$ of
  $\Gamma$ there exists a directed path from $u_1$ to $u_2$ in
  $\Gamma$).
\end{enumerate}

\end{lem}

We can now give a precise description of those rational points in
$Q_m$ that come from some actual cyclic words in $F$.

\begin{thm}\label{realize}
  Let $m\ge 2$ and $q\in Q_m$.  Then $q\in \alpha_m(C)$ if and only if
  the underlying topological graph of $\Gamma_q'$ is connected.
\end{thm}

\begin{proof}
  Suppose $q\in Q_m$ is a rational point where the underlying graph of
  $\Gamma_q'$ is connected.
  
  Choose a positive integer $N>0$ such that for every $v\in F$ with
  $|v|=m$ we have $Nq_v\in \mathbb Z$.

  Let $\Delta:=[N\Gamma_q']$. Thus the vertex sets of $\Delta$ and of
  $\Gamma_q'$ are the same.
  
  Then $\Delta$ is a connected digraph with $N$ directed edges where
  for each vertex the in-degree is equal to the out-degree.  Therefore
  $\Delta$ possesses an Euler circuit (in the di-graph sense).
  
  Choose such an Euler-circuit and also, by specifying an initial
  vertex, choose a closed edge-path $e_1,\dots, e_N$ in $\Delta$
  realizing this circuit.  Each $e_i$ corresponds to a freely reduced
  word $v_i$ of length $m$ such that the origin of $e_i$ is
  $(v_i)_{-}$ and the terminus of $e_i$ is $(v_i)_+$.
  
  We will create a cyclically reduced word $z$ of length $N$ as
  follows. Put $z_1$ to be the last letter of $v_1$. If the words
  $z_1,\dots, z_{i-1}$ are already constructed, then we define $z_i$
  to be the word $z_ix$ where $x$ is the last letter of $v_i$ (and so
  of $(v_i)_+$). Put $z:=z_t$. It is easy to see that $z$ is
  cyclically reduced. Let $w=(z)$ be the cyclic word defined by $z$.
  Then by construction we have $||w||=N$. Moreover, it is easy to see
  that for every $v\in F$ with $|v|=m$ we have $n_w(v)=Nq_v$,
  $f_w(v)=q_v$. Hence $\alpha_m(w)=q$.
  
  Note that in the above construction the cyclic word $w$ depends only
  on the Euler circuit in $\Delta$ but not on an actual
  closed-edge-path realizing this circuit.

  Suppose now that $q=\alpha_m(w)\in Q_m$ for some nontrivial cyclic
  word $w$ in $F$. Put $N=||w||$. We need to show that the graph
  $\Gamma_q'$ is connected.

  For each $v\in F$ with $|v|=m$ we have $n_w(v)=Nw_v'\in \mathbb Z$.
  Hence with this choice of $N$ the graph $N\Gamma_q'$ has positive
  integer edge-labels whose sum is equal to $N$.  We will produce a
  path in $N\Gamma_q$ passing through each edge of $N\Gamma_q'$.  This
  would imply that $\Gamma_q'$ is connected.

  Recall that the edge-set of $\Gamma_q'$ (and hence of $N\Gamma_q'$)
  is in a canonical one-to-one correspondence with the set $\{v\in F :
  |v|=m, f_w(v)\ne 0\}$. Choose a cyclically reduced word $z$ with
  $(z)=w$. Construct an edge-path $e_1,\dots, e_N$ in $N\Gamma_q'$ as
  follows.
  
  Write $z$ as $z=x_1 x_2\dots x_{N}$ where $x_j\in X$. Also, for
  $t\ge N$ put $x_t=x_i$ where $1\le i\le N$ and $t\equiv i \mod N$.
  For $i=1,\dots, N$ put $e_i$ to be the edge in $N\Gamma_q'$
  corresponding to $x_i\dots x_{i+m-1}$, that is an edge from
  $x_i\dots x_{i+m-2}$ to $x_{i+1}\dots x_{i+m-1}$. It is easy to see
  that $e_1,\dots, e_N$ is a closed edge-path in $N\Gamma_q'$, where
  each edge corresponding to a word $v\in F$ with $|v|=m$ is repeated
  exactly $Nf_w(v)=n_w(v)$ times. Since this path passes through every
  edge of $N\Gamma_q$, the underlying graph of $\Gamma_q'$ is
  connected, as claimed.
\end{proof}

\begin{exmp}
  Let $k=2$, $F=F(a,b)$ and $m=3$.
  
  Consider a point $q\in Q_3$ given by $q_{aba}=2/5$,
  $q_{bab}=q_{a^2b}=q_{ba^2}=1/5$ and all other $q_v=0$, $|v|=3$,
  $v\not\in\{aba, bab, a^2b, ba^2 \}$.
  
  Take $N=5$. Then the graph $5\Gamma_q'$ has the vertices $ab, ba,
  a^2$ and four directed edges: the edge $e_1$ from $ab$ to $ba$ with
  label $2$, the edge $e_2$ from $ba$ to $ab$ with label $1$, the edge
  $e_3$ from $ba$ to $a^2$ with label $1$ and the edge $e_4$ from
  $a^2$ to $ab$ with label $1$. Then the path $e_1,e_2,e_1,e_3,e_4$
  defines an Euler cycle at a vertex $ab$ in the graph $[5\Gamma_q']$.
  The corresponding cyclically reduced word $z$ is
\[
z=abaab
\]
Then for $w=(z)$ we have $||w||=5$, $n_w(aba)=2$,
$n_w(bab)=n_w(ba^2)=n_w(a^2b)=1$ and $\alpha_3(w)=q$.

Similarly, the path
\[
e_1,e_2,e_1,e_2,e_1,e_3,e_4,e_1,e_3,e_4
\]
defines an Euler cycle at vertex $ab$ in the graph $[10\Gamma_q']$.
The corresponding cyclically reduced word $z_1$ is
\[
z_1=ababaabaab
\]
Then for $w_1=(z_1)$ we have $||w||=10$, $n_w(aba)=4$,
$n_w(bab)=n_w(ba^2)=n_w(a^2b)=2$ and $\alpha_3(w)=q$.
\end{exmp}

\begin{lem}\label{real1}
  Let $q\in Q_1$ be a rational point. Then $q\in \alpha_1(C)$ if and
  only if there does not exist a letter $x\in X$ such that $q_x>0$,
  $q_{x^{-1}}>0$ and $q_x+q_{x^{-1}}=1$ (and hence $q_y=0$ for all
  $y\in X$ such that $y\ne x$, $y\ne x^{-1}$).

\end{lem}

\begin{proof}
  
  To simplify notation we denote $A_i=a_i^{-1}$ throughout the proof.
  It is clear that if $q=\alpha_1(w)$ for a nontrivial cyclic word $w$
  then there does not exist a letter $x\in X$ such that $q_x>0$,
  $q_{x^{-1}}>0$ and $q_x+q_{x^{-1}}=1$.

  Suppose now that $q\in Q_1$ is a rational point such that there does
  not exist a letter $x\in X$ with $q_x>0$, $q_{x^{-1}}>0$ and
  $q_x+q_{x^{-1}}=1$. We need to show that $q$ can be realized as
  $\alpha_1(w)$ for some cyclic word $w$.  Choose an integer $N>0$
  such that all coordinates of the point $Nq$ are integers, that is
  $Nq_a\in \mathbb Z$ for each $a\in \{a_1,\dots a_k\}^{\pm 1}$.

  We define the words $\alpha_i, \beta_i$ for $i=1,\dots, k$ as
  follows.

\begin{itemize}
\item If $q_{a_i}=q_{A_i}=0$ we put $\alpha_i=\beta_i=1$. \item If
  $q_{a_i}>0, q_{A_i}=0$, we put $\alpha_i=\beta_i=a_i^{Nq_{a_i}}$.
\item If $q_{A_i}>0, q_{a_i}=0$, we put
  $\alpha_i=\beta_i=A_i^{Nq_{a_i}}$.
\item If $q_{a_i}>0, q_{A_i}>0$, we put $\alpha_i=a_i^{2Nq_{a_i}}$,
  $\beta_i=A_i^{2Nq_{A_i}}$.
\end{itemize}

Now put $v=\alpha_1\dots\alpha_k\beta_1\dots \beta_k$. It is clear
that $v$ is a cyclically reduced word of length $2N$. For the cyclic
word $w=(v)$ we have $n_w(a)=2Nq_a$, $f_w(a)=q_a$ for every letter $a$
and hence $q=\alpha_1(w)\in \alpha_1(C)$, as required.

\end{proof}

\begin{conv}
  For $m\ge 1$ let $Q_m^+$ be the set of all points $q\in Q_m$ such
  that for every $v\in F$ with $|v|=m$ we have $q_v>0$. Let $QQ_m^+$
  be the set of rational points in $Q_m^+$.
\end{conv}

It is clear that if $m\ge 1$ then for any $q\in Q_m^+$ the graph
$\Gamma_1'$ is connected and hence by Theorem~\ref{realize} and
Lemma~\ref{real1}, $QQ_m^+\subseteq Q_m^+\subseteq \alpha_m(C)$.

\begin{lem}\label{dense}
  Let $m\ge 1$. Then the set $QQ_m^+$ is dense in $Q_m$.
\end{lem}
\begin{proof}
  
  The statement is obvious when $m=1$, so we will assume that $m\ge
  2$.
  
  Since $Q_m$ is a compact convex finite-dimentional polyhedron defined by equations and inequalities with rational coefficients, the set of rational points $QQ_m$ is dense in $Q_m$. Thus it suffices to show that $QQ_m^{+}$ is dense in $QQ_m$.

  Recall that $D(m)$ is the number of freely reduced words of length $m$ in $F$. Put $z_u=\frac{1}{D(m)}$ for all $u\in F$ with $|u|=m$ and put $z=(z_u)_{|u|=m}$. Then it is easy to see that $z\in Q_m$ and in fact $z\in QQ_m^+$.

  Let $q\in QQ_m$ be arbitrary. Let $1>\epsilon>0$ be a rational number. By convexity of $Q_m$ we have $\epsilon z +(1-\epsilon) q\in Q_m$ and, moreover $\epsilon z +(1-\epsilon) q\in QQ_m^+$. Since
\[
\lim_{\epsilon\to 0} \epsilon z +(1-\epsilon) q=q,
\] 
and $q\in QQ_m$ was arbitrary, it follows that $QQ_m^+$ is dense in $QQ_m$ as required.
\end{proof}

Note that Lemma~\ref{dense} and the definition of $Q(F)$ immediately
imply that $\alpha(C)$ is dense in $Q(F)$, yielding another proof of
Corollary~\ref{cor:dense}

\begin{cor}\label{onto}
  If $m\ge 1$ then $\pi_m(Q_m)=Q_{m-1}$.
\end{cor}

\begin{proof}
  The case $m=1$ is obvious, so we assume $m\ge 2$.  It follows from
  the definition of $\pi_m$ that $\pi_m(Q_m^+)\subseteq Q_{m-1}^+$ and
  $\pi_m(QQ_m^+)\subseteq QQ_{m-1}^+$.
  
  Moreover, since $Q_m$ and $Q_{m-1}$ are compact and $QQ_m^+$ is
  dense in $Q_m$, the image of $\pi_m$ is equal to the closure of
  $\pi_m(QQ_m^+)$. Suppose $\pi_m$ is not onto, so that $\pi_m(Q_m)$
  is a proper compact convex sub-polyhedron is a compact convex
  polyhedron $Q_{m-1}$. Therefore $Q_{m-1}-\pi_m(Q_m)$ contains a
  nonempty open subset of $Q_{m-1}$ and therefore, by
  Lemma~\ref{dense}, it contains a point $q$ of $QQ_{m-1}^+$.
  
  By Theorem~\ref{realize} and Lemma~\ref{real1} there exists a cyclic
  word $w$ such that $\alpha_{m-1}(w)=q$. Then $q=\pi_m(\alpha_m(w))$
  and hence $q\in \pi_m(Q_m)$, contrary to our assumption that $q\in
  Q_{m-1}-\pi_m(Q_m)$.
\end{proof}

\begin{lem}\label{ext}
Let $m\ge 2$ and let $p\in Q_m$ be an extremal point of the polyhedron
$Q_m$. Then $\Gamma'_p$ is connected and hence, by
Theorem~\ref{realize}, $p\in \alpha_m(C)$.
\end{lem}
\begin{proof}
Note that all $p_u$ (where $u\in F, |u|=m$) are rational numbers since
$p$ is an extremal point of a compact comvex polyhedron $Q_m$ defined
by equations and inequalities with integer coefficients.

Suppose that $\Gamma_p'$ is not connected.
Then there exists two nonempty disjoint labelled subgraphs $G_1,G_2$ of $\Gamma_p'$ such
that $\Gamma_p'=G_1\sqcup G_2$.
Recall that $G_p'$ consists of all the edges of $G_p$ with positive
labels and of their endpoints. Let $s$ be the sum of the edge-labels
in $G_1$ and let $t$ be the sum of the edge-labels in $G_2$. Then
$s+t=1, s>0, t>0$.

Define $r=(r_u)_{|u|=m}$ as follows.
Put $r_u=0$ if $u\in F, |u|=m$ does not correspond to an edge of
$G_1$. Put $r_u=\frac{p_u}{s}$ if $u\in F, |u|=m$ does correspond
to an edge of $G_1$. Then $r\in Q_m$ and $\Gamma_r'=G_1$.

Similarly, define  $q=(q_u)_{|u|=m}$ as follows.
Put $q_u=0$ if $u\in F, |u|=m$ does not correspond to an edge of
$G_2$. Put $q_u=\frac{p_u}{t}$ if $u\in F, |u|=m$ does correspond
to an edge of $G_2$. Then $q\in Q_m$ and $\Gamma_q'=G_2$.

By construction we have $sr+tq=p$ and $r\ne p, q\ne p$. Thus $q\in
Q_m$ is a convex linear combination of two points of $Q_m$ different
from $p$. This contradicts our assumption that $p$ is an extremal
point of $Q_m$.
\end{proof}

In order to estimate the topological dimension of $Q_m$ we need the
following elementary fact:

\begin{lem}\label{lem:dim}
  Let $f_1,\dots, f_s, g_1\dots, g_t:\mathbb R^n\to \mathbb R$ be
  affine functions.  Let \[J=\{x\in {\mathbb R^n}: f_i(x)=0, 1\le i\le
  s\}\] and \[Z=\{x\in {\mathbb R^n}: f_i(x)=0, g_j(x)>0, 1\le i\le s,
  1\le j\le t\}.\] Suppose $Z$ is nonempty. Then ${\rm dim}\ Z={\rm
    dim}\ J$.
\end{lem}

\begin{cor}
  Let $m\ge 2$. Then
\[
D(m)-D(m-1)-1\le {\rm dim}\ Q_m \le D(m)-1.
\]
\end{cor}
\begin{proof}
  The upper bound is obvious from condition (1) of
  Definition~\ref{defn:qm}. Let $J_m$ be the subset of ${\mathbb
    R}^{D(m)}$ defined by equations from parts (1) and (3) in
  Definition~\ref{defn:qm}. Thus $J_m$ is defined by $D(m-1)+1$ affine
  equations, of which $D(m-1)$ equations (those coming from part (3))
  are homogeneous. Let $J_m'$ be the set defined by equations from
  part (3) of Definition~\ref{defn:qm} and the homogeneous version of
  the equation from part (1), namely $\sum_{|v|=m} q_v=0$. The set
  $J_m'$ is given by $D(m-1)+1$ linear homogeneous equations in $D(m)$
  variables. Hence ${\rm dim}\ J_m'\ge D(m)-D(m-1)-1$. Since $Q_m$ is
  nonempty, we can choose $q_0\in Q_m\subseteq J_m$. Then
  $J_m=q_0+J_m'$ and hence ${\rm dim}\ J_m={\rm dim}\ J_m'\ge
  D(m)-D(m-1)-1$. Lemma~\ref{lem:dim} implies that ${\rm dim}\ 
  J_m={\rm dim}\ Q_m$, so that ${\rm dim}\ Q_m\ge D(m)-D(m-1)-1$, as
  required.
\end{proof}

Note that directly from Definition~\ref{defn:qm} we get ${\rm dim}\ 
Q_1=D(1)-1=2k-1$.

We summarize our results regarding $Q_{m}$ and $Q(F)$ accumulated so
far:

\begin{thm} The following hold:
\begin{enumerate}
\item If $m\ge 1$ then $\pi_m(Q_m)=Q_{m-1}$.
  
\item We have $QQ_m^+\subseteq \alpha_m(C)\subseteq Q_m$ and $QQ_m^+$
  is dense in $Q_m$ for $m\ge 0$
  
\item The space $Q(F)$ with the inverse limit topology is compact.
  
\item The set $\alpha(C)$ is dense in $Q(F)$.
  
\item If $w\in C$ is a nontrivial cyclic word then for any $n\ge 1$ we
  have $\alpha(w)=\alpha(w^n)$ and $\alpha_m(w)=\alpha_m(w^n)$, where
  $m\ge 1$.
  
\item Let $w, u$ be nontrivial cyclic words. Then
  $\alpha(w)=\alpha(u)$ if and only if both $w$ and $u$ are positive
  powers of the same cyclic word.
  
\item For $q\in Q_m$ is a rational point, where $m\ge 2$, then we have
  $q\in \alpha_m(C)$ if and only if the underlying undirected graph of
  $\Gamma_q'$ is connected.
  
\item If $q\in Q_1$ is a rational point then $q\in \alpha_1(C)$ if and
  only if there does not exist $x\in X$ such that $q_x+q_{x^{-1}}=1$,
  $q_x>0$, $q_{x^{-1}}>0$.
  
\item We have
\[
D(m)-D(m-1)-1\le {\rm dim}\ Q_m \le D(m)-1
\]
for any $m\ge 2$. Also we have ${\rm dim}\ Q_1=D(1)-1=2k-1$ and ${\rm
  dim}\ Q_0=0$
\item If $m\ge 2$ and $p\in Q_m$ is an extremal point of $Q_m$ then
  $p\in \alpha_m(C)$.
\end{enumerate}
\end{thm}

\section{The action of $Out(F)$ on the frequency space}

\begin{defn}[Nielsen automorphisms]
  A \emph{Nielsen} automorphism of $F$ is an automorphism $\tau$ of
  one of the following types:

\begin{enumerate}
  
\item There is some $i, 1\le i\le k$ such that $\tau(a_i)=a_i^{-1}$
  and $\tau(a_j)=a_j$ for all $j\ne i$.
  
\item There are some $1\le i<j\le k$ such that $\tau(a_i)=a_j$,
  $\tau(a_j)=a_i$ and $\tau(a_l)=a_l$ when $l\ne i, l\ne j$.
  
\item There are some $1\le i<j\le k$ such that $\tau(a_i)=a_ia_j$ and
  $\tau(a_l)=a_l$ for $l\ne i$.

\end{enumerate}

\end{defn}
It is well-known~(see, for example, \cite{LS}) that $Aut(F)$ is
generated by the set of Nielsen automorphisms.  The following two
propositions are crucial for our arguments (as well as for our results
in \cite{KSS}).

\begin{prop}\label{nil}
  For each Nielsen automorphism $\tau$ of $F$ and for any $u\in F$
  with $|u|=m$ there exist integers $c(v)=c(v,u,\tau)\ge 0$, where
  $v\in F, |v|=2m+6$, such that for any cyclic word $w$ we have
\[
n_{\tau(w)}(u)=\sum_{|v|=2m+6} c(v) n_w(v).
\]
\end{prop}

\begin{proof}
  
  Note that in view of part (1) of Lemma~\ref{red} it suffices to
  express $n_{\tau(w)}(u)$ as a linear function in terms of $n_w(v)$
  where $|v|\le 2m+6$.

  The statement of the proposition is obvious if $\tau$ interchanges
  two generators or inverts a generator since in that case
  $n_{\tau(w)}(u)=n_w(\tau^{-1}(u))$.
  
  Without loss of generality we may now assume that
  $F=F(a,b,a_3,...,a_k)$ with $a=a_1, b=a_2$ and that $\tau(a)=a,
  \tau(b)=ba$, $\tau(a_i)=a_i$ for $2<i\le k$.
  
  Note that if $i\ne 1$ then no occurrences of $a_i^{\pm 1}$ cancel
  when we reduce $\tau(w)$. The cancellation structure involving
  $a^{\pm 1}$ is also very controlled.  Only the occurrences of
  $ba^{-1}$ and $ab^{-1}$ in $w$ produce appearances of $aa^{-1}$ in
  the non-reduced form of $\tau(w)$. There are no appearances of
  $a^{-1}a$ in the non-reduced form of $\tau(w)$. We will say that an
  occurrence of $a$ in the reduced form of $\tau(w)$ is \emph{old} if
  it comes from an occurrence of $a$ in $w$ that was not followed by
  $b^{-1}$ in $w$ (and hence was not cancelled in the process of
  freely reducing $\tau(w)$.) We will say that an occurrence of $a$ in
  the reduced form of $\tau(w)$ is \emph{new} if it comes from an
  occurrence of $b$ in $w$ that was not followed by $a^{-1}$ in $w$
  (and hence was not cancelled in the process of reducing $\tau(w)$.)
  Similarly, an occurrence of $a^{-1}$ in the reduced form of
  $\tau(w)$ is \emph{old} if it comes from an occurrence of $a^{-1}$
  in $w$ that was not preceded by $b$ in $w$.  An occurrence of
  $a^{-1}$ in the reduced form of $\tau(w)$ is \emph{new} if it comes
  from an occurrence of $b^{-1}$ in $w$ that was not preceded by $a$
  in $w$.

  The argument involves a straightforward (but tedious) analysis of
  several cases for the possible configurations at the beginning and
  the end of $u$. We present the details for completeness.
  In various summations below we will implicitly assume that $x$
  varies over the set $\{a,b, a_3,\dots, a_k\}^{\pm 1}$. We will also
  implicitly assume that the summation occurs only over those words
  satisfying the explicitly specified conditions that are freely
  reduced.

\noindent{\bf Case 1.} Let $u=a^s$ where $s>0$.

The number of those occurrences of $a^s$ in $\tau(w)$ that start with
an old occurrence of $a$ is equal to $\sum_{x\ne b^{-1}}n_w(a^sx)$.
The number of those occurrences of $a^s$ in $\tau(w)$ that start with
a new occurrence of $a$ is equal to $\sum_{x\ne
  b^{-1}}n_w(ba^{s-1}x)$.

Thus

\[
n_{\tau(w)}(u)=\sum_{x\ne b^{-1}}n_w(a^sx)+\sum_{x\ne
  b^{-1}}n_w(ba^{s-1}x).
\]

\noindent{\bf Case 2.} Let $u=a^{-s}$ where $s>0$.

The number of those occurrences of $a^{-s}$ in $\tau(w)$ that start
with an old occurrence of $a^{-1}$ is equal to $\sum_{x\ne
  b}n_w(xa^{-s})$. The number of those occurrences of $a^{-s}$ in
$\tau(w)$ that start with a new occurrence of $a^{-1}$ is equal to $0$
if $s>1$ and is equal to $\sum_{x\ne a}n_w(xb^{-1})$ if $s=1$.

Thus for $s=1$
\[
n_{\tau(w)}(a^{-1})=\sum_{x\ne b}n_w(xa^{-1})+\sum_{x\ne
  a}n_w(xb^{-1})
\]

and for $s>1$
\[
n_{\tau(w)}(a^{-s})=\sum_{x\ne b}n_w(xa^{-s}).
\]

\noindent{\bf Case 3.} Let $u=a^sza^t$ where $s\in \mathbb Z$,
$t\in \mathbb Z$ and $z$ is a nontrivial reduced word that neither
begins nor ends with $a^{\pm 1}$.

Let $z'$ be obtained from the reduced form of $\tau^{-1}(z)$ by
deleting the first letter if that letter is $a^{\pm 1}$ and by
deleting the last letter if that letter is $a^{\pm 1}$. Note that now
the first letters of $z'$ and $z$ are the same and the last letters of
$z'$ and $z$ are the same.

\noindent{\bf Subcase 3.1} Suppose that $s=t=0$, so that $u=z$.

It is easy to see that $n_{\tau(w)}(u)=n_w(z')$.

\noindent{\bf Subcase 3.2} Suppose that $s>0, t>0$.

\noindent {\bf Subcase 3.2.A} Assume that $z$ does not end in $b$.

If $z$ does not start with $b^{-1}$ then the number of occurrences of
$u$ in $w$ that start with an old occurrence of $a$ is equal to
$\sum_{x\ne b^{-1}} n_w(a^sz'a^tx)$. In this case the number of
occurrences of $u$ in $w$ that start with a new occurrence of $a$ is
equal to $\sum_{x\ne b^{-1}} n_w(ba^{s-1}z'a^tx)$. Thus if $z$ does
not start with $b^{-1}$ then

\[
n_{\tau(w)}(u)=\sum_{x\ne b^{-1}} n_w(a^sz'a^tx)+\sum_{x\ne b^{-1}}
n_w(ba^{s-1}z'a^tx).
\]

Assume now that $z$ starts with $b^{-1}$. Then the number of
occurrences of $u$ in $w$ that start with an old occurrence of $a$ is
equal to $\sum_{x\ne b^{-1}} n_w(a^{s+1}z'a^tx)$. In this case the
number of occurrences of $u$ in $w$ that start with a new occurrence
of $a$ is equal to $\sum_{x\ne b^{-1}} n_w(ba^{s}z'a^tx)$. Thus if $z$
starts with $b^{-1}$ then

\[
n_{\tau(w)}(u)=\sum_{x\ne b^{-1}} n_w(a^{s+1}z'a^tx)+\sum_{x\ne
  b^{-1}} n_w(ba^{s}z'a^tx).
\]

{\bf Subcase 3.2.B} Suppose now that $z$ ends in $b$.

Assume first that $z$ does not begin with $b^{-1}$.

Then the number of occurrences of $u$ in $w$ that start with an old
occurrence of $a$ is equal to $\sum_{x\ne b^{-1}} n_w(a^sz'a^{t-1}x)$.
In this case the number of occurrences of $u$ in $w$ that start with a
new occurrence of $a$ is equal to $\sum_{x\ne b^{-1}}
n_w(ba^{s-1}z'a^{t-1}x)$. Thus if $z$ does not start with $b^{-1}$
then

\[
n_{\tau(w)}(u)=\sum_{x\ne b^{-1}} n_w(a^sz'a^{t-1}x)+\sum_{x\ne
  b^{-1}} n_w(ba^{s-1}z'a^{t-1}x).
\]

Assume now that $z$ does begin with $b^{-1}$. Then the number of those
occurrences of $u$ in $w$ that start with an old occurrence of $a$ is
equal to $\sum_{x\ne b^{-1}} n_w(a^{s+1}z'a^{t-1}x)$. In this case the
number of occurrences of $u$ in $w$ that start with a new occurrence
of $a$ is equal to $\sum_{x\ne b^{-1}} n_w(ba^{s}z'a^{t-1}x)$. Thus if
$z$ starts with $b^{-1}$ then

\[
n_{\tau(w)}(u)=\sum_{x\ne b^{-1}} n_w(a^{s+1}z'a^{t-1}x)+\sum_{x\ne
  b^{-1}} n_w(ba^{s}z'a^{t-1}x).
\]

\noindent{\bf Subcase 3.3.} Suppose that $s>0$ and $t=0$ so that
$u=a^sz$.

Then, as in Subcase 3.2, if $z$ does not begin with $b^{-1}$ then

\[
n_{\tau(w)}(u)=n_w(a^sz')+n_w(ba^{s-1}z').
\]

Similarly, if $z$ begins with $b^{-1}$ then arguing as in Subcase 3.2
we get

\[
n_{\tau(w)}(u)= n_w(a^{s+1}z')+n_w(ba^{s}z').
\]

\noindent{\bf Subcase 3.4.} Suppose that $s=0$ and $t>0$ so that
$u=za^t$.

It is easy to see that if $z$ does not ends in $b$ then

\[
n_{\tau(w)}(u)=\sum_{x\ne b^{-1}} n_w(z'a^tx).
\]

Similarly, if $z$ ends in $b$ then
\[
n_{\tau(w)}(u)=\sum_{x\ne b^{-1}} n_w(z'a^{t-1}x).
\]

\noindent{\bf Subcase 3.5} We will consider in detail one more
case and leave the remaining possibilities to the reader.

Suppose that $u=a^{-s}za^{-t}$ where $s,t>0$ and $z$ is a nontrivial
reduced word that neither begins nor ends with $a^{\pm 1}$.

{\bf Subcase 3.5.A} Suppose that $z$ does not end in $b$.

If $z$ does not start with $b^{-1}$ then

\[
n_{\tau(w)}(u)=\sum_{x\ne b} n_w(xa^{-s}z'a^{-t}).
\]

If $z$ starts with $b^{-1}$ and $s>1$ then

\[
n_{\tau(w)}(u)=\sum_{x\ne b} n_w(xa^{-s+1}z'a^{-t}).
\]

If $z$ starts with $b^{-1}$ and $s=1$, so that $u=a^{-1}za^{-t}$, then

\[
n_{\tau(w)}(u)=\sum_{x\ne a} n_w(xz'a^{-t}).
\]

{\bf Subcase 3.5.B} Suppose that $z$ ends in $b$.

If $z$ does not start with $b^{-1}$ then

\[
n_{\tau(w)}(u)=\sum_{x\ne b} n_w(xa^{-s}z'a^{-t-1}).
\]

If $z$ starts with $b^{-1}$ and $s>1$ then

\[
n_{\tau(w)}(u)=\sum_{x\ne b} n_w(xa^{-s+1}z'a^{-t-1}).
\]

If $z$ starts with $b^{-1}$ and $s=1$ then

\[
n_{\tau(w)}(u)=\sum_{x\ne a} n_w(xz'a^{-t-1}).
\]

\end{proof}

\begin{prop}\label{phi}
  Let $\phi\in Out(F)$ be an outer automorphism and let $t$ be such
  that $\phi$ can be represented, modulo $Inn(F)$, as a product of $t$
  Nielsen automorphisms.
  
  Then for any freely reduced word $v\in F$ with $|v|=m$ there exists
  a collection of integers $c(u,v)=c(u,v,\phi)\ge 0$, where $u\in F$,
  $|u|=8^tm$, such that for any nontrivial cyclic word $w$ we have

\[
n_{\phi(w)}(v)=\sum_{|u|=8^tm} c(u,v) n_w(u).
\]
\end{prop}

\begin{proof}
  Let $\phi=\tau_1\dots \tau_t\in Out(F)$ where each $\tau_i$ is a
  Nielsen automorphism. We prove Proposition~\ref{phi} by induction on
  $t$.
  
  For $t=1$ Proposition~\ref{phi} follows from Proposition~\ref{nil}
  since $2m+6\le 8m$ for $m\ge 1$. Suppose now that $t>1$ and that
  Proposition~\ref{phi} has been established for $t-1$.
  
  Put $\psi= \tau_2 \dots \tau_t$.
  
  By assumption there is a collection of numbers $b(v,z)\ge 0$, where
  $|v|=8^{t-1}(2m+6), |z|=2m+6$, such that for any cyclic word $w$ we
  have

\[
n_{\psi(w)}(z)=\sum_{|v|= 8^{t-1}(2m+6)} b(v,z) n_w(v).
\]

Since $\phi(w)=\tau_1(\psi(w))$, Proposition~\ref{nil} implies that

\begin{align*}
  n_{\phi(w)}(u)&=\sum_{|z|=2m+6} c(z,u,m) n_{\psi(w)}(z)=\\
  &=\sum_{|z|=2m+6} c(z,u,m) \sum_{|v|=8^{t-1}(2m+6)} b(v,z)
  n_w(v)=\\
  &=\sum_{|v|=8^{t-1}(2m+6)} [\sum_{|z|= 2m+6} b(v,z) c(z,u,m)]
  n_w(v).
\end{align*}
Since $8^{t-1}(2m+6)\le 8^tm$, in view of part (1) of Lemma~\ref{red},
the statement of Proposition~\ref{phi} follows.
\end{proof}

\begin{prop}\label{phi1}
  Let $\phi\in Out(F)$ be an outer automorphism and let $t$ be such
  that $\phi$ can be represented, modulo $Inn(F)$, as a product of $t$
  Nielsen automorphisms.
  
  Then there exists
  a collection of integers $d(u,\phi)\ge 0$, where $u\in F$,
  $|u|=8^t$, such that for any nontrivial cyclic word $w$ we have
\[
f_{\phi(w)}(v)=\frac{\sum_{|u|=8^tm} c(u,v) f_w(u)}{\sum_{|u|=8^t}
  d(u,\phi) f_w(u)}
\]
and
\[
||\phi(w)||=\sum_{|u|=8^t} d(u,\phi) n_w(u)
\]

Here $c(u,v)=c(u,v,\phi)$ are as in Proposition~\ref{phi}.

Moreover, there is $c_0>0$ such that for any $q\in Q_{8^t}$
\[
\sum_{|u|=8^t} d(u,\phi) q_u\ge c_0>0.
\]

\end{prop}
\begin{proof}
  By Proposition~\ref{phi} we have
\[
n_{\phi(w)}(v)=\sum_{|u|=8^tm} c(u,v) n_w(u)
\]
and

\[
||\phi(w)||=\sum_{|z|=1} n_{\phi(w)}(z)=\sum_{|z|=1}\sum_{|u|=8^t}
c(u,z) n_w(u).
\]
Therefore

\begin{gather*}
  f_{\phi(w)}(v)=\frac{n_{\phi(w)}(v)}{||\phi(w)||}=\frac{\sum_{|u|=8^tm}
    c(u,v) n_w(u)}{\sum_{|z|=1}\sum_{|u|=8^t} c(u,z) n_w(u)}=\\
  \frac{\sum_{|u|=8^tm} c(u,v)
    n_w(u)/||w||}{\sum_{|z|=1}\sum_{|u|=8^t} c(u,z)
    n_w(u)/||w||}=\frac{\sum_{|u|=8^tm} c(u,v) f_w(u)}{\sum_{|u|=8^t}
    [\sum_{|z|=1} c(u,z)] f_w(u)},
\end{gather*}
and the first part of the statement follows with $d(u)=\sum_{|z|=1}
c(u,z)$.

Note that by construction

\[
\frac{||\phi(w)||}{||w||}=\sum_{|u|=8^t}
    [\sum_{|z|=1} c(u,z)] f_w(u)=\sum_{|u|=8^t} d(u) f_w(u),
\]
for any $w\in C$. Choose an automorphism $\psi\in Aut(F)$ representing
$\phi\in Out(F)$ (so that $\psi^{-1}$ represents $\phi^{-1}$).

Obviously $\frac{||\phi(w)||}{||w||}\le L(\psi)$ where $L(\psi)=\max
\{|\psi(x)|: x\in X\}$ and $1\le L(\psi)<\infty$.  Therefore
$\frac{||\phi(w)||}{||w||}\ge \frac{1}{L(\psi^{-1})}$ for any $w\in
C$. Hence the function $h(q)=\sum_{|u|=8^t} d(u) q_u$ on $Q_{8^t}$
satisfies $h(q)\ge \frac{1}{L(\psi^{-1})}$ for each $q\in
\alpha_{8^t}(C)$. Since $\alpha_{8^t}(C)$ is dense in $Q_{8^t}$,
this implies that $h(q)\ge \frac{1}{L(\psi^{-1})}>0$ for every $q\in
Q_{8^t}$, as claimed.
\end{proof}

\begin{conv}
  For each $\phi\in Out(F)$ we choose and fix a shortest
  representation of $\phi$ modulo $Inn F$ as a product Nielsen
  automorphisms. Denote by $t=t(\phi)$ the number of automorphisms in
  this product.  We also fix the numbers $c(u,v,\phi)$ (where $v,u,\in
  F$, $|u|=8^t|v|$) provided by Proposition~\ref{phi}. For any cyclic
  word $w$ we will denote $\alpha(w)$ by $w^{\alpha}$.
\end{conv}

\begin{thm}\label{action}
  For $\phi\in Out(F)$ and $q\in \alpha(C)\subseteq Q(F)$ put
  $\tilde\phi(q):=\phi(w)^{\alpha}$, where $w\in C$ is any cyclic word
  with $w^{\alpha}=q$.
  
  Then $\tilde\phi(q)$ is well-defined and the map $\tilde\phi:
  \alpha(C)\to \alpha(C)$ is continuous in the induced from $Q(F)$
  topology. This map extends uniquely to a continuous homeomorphism
  $\tilde\phi: Q(F)\to Q(F)$.  The map
\[\tilde{}:
Out(F)\to Homeo(Q(F)), \quad \tilde{ }: \phi\mapsto \tilde\phi, \text{
  for each } \phi\in Out(F),\] defines a left action of $Out(F)$ on
$Q(F)$ by homeomorphisms.
\end{thm}

\begin{proof}
  By Proposition~\ref{imalpha} we know that if
  $q=w_1^{\alpha}=w_2^{\alpha}$ then $w_1$, $w_2$ are positive powers
  of the same cyclic word and therefore the same is true for
  $\phi(w_1)$ and $\phi(w_2)$. Hence
  $\phi(w_1)^{\alpha}=\phi(w_2)^{\alpha}$ and so $\tilde\phi:
  \alpha(C)\to \alpha(C)$ is well-defined.
  
  To see that this map is continuous and extends uniquely to a
  continuous map $Q(F)\to Q(F)$, choose an arbitrary integer $m>0$.
  Then for any $q=w^{\alpha}$ and any $v\in F$ with $|v|=m$ we have

\[
\tilde\phi(q)_v=f_{\phi(w)}(v)=\frac{\sum_{|u|=8^tm} c(u,v,\phi)
  f_{w}(u)}{\sum_{|u|=8^t} d(u,\phi) f_{w}(u)}
=\frac{\sum_{|u|=8^tm} c(u,v,\phi) q_u}{\sum_{|u|=8^t} d(u,\phi)
  q_u}.
\]

Fix $v\in F$ with $|v|=m$. Then by Proposition~\ref{phi1} the
fractional-linear function
\[
r_v(q):=\frac{\sum_{|u|=8^tm} c(u,v,\phi) q_u}{\sum_{|u|=8^t}
  d(u,\phi) q_u}\
\]
has denominator satisfying
\[
\sum_{|u|=8^t} d(u,\phi) q_u\ge c_0>0
\]
on $Q(F)$ for some constant $c_0>0$ independent of $q\in
Q(F)$. Hence $r_v(q)$ is continuous on $Q(F)$.

Since $\alpha(C)$ is dense in $Q(F)$, this implies that the above map
$\tilde\phi: \alpha(C)\to \alpha(C)$ extends uniquely to a continuous
map $\tilde\phi:Q(F)\to Q(F)$.  The map $\tilde{}: \phi\mapsto
\tilde\phi$ defines an action of $Out(F)$ on $Q(F)$ by homeomorphisms.
Indeed, since $(\phi\psi)(w)=\phi(\psi(w))$ for any $w\in C$ and
$\phi,\psi\in Out(F)$, we conclude that $\tilde \phi \circ \tilde
\phi= \tilde(\phi\psi)$ on $\alpha(C)$ and therefore, by continuity,
on $Q(F)$. Also, for the trivial automorphism $1\in Out(F)$ the action
of $1$ on $C$ is trivial, and hence $\tilde
1|_{\alpha(C)}=Id_{\alpha(C)}$. Therefore, by continuity $\tilde
1|_{Q(F)}=Id_{Q(F)}$. Thus $\tilde{ }: \phi:\to \tilde \phi$, $\phi\in
Out(F)$ indeed is an action of $Out(F)$ on $Q(F)$ by homeomorphisms.
\end{proof}

\begin{rem}
As we mentioned in the introduction, the formulas defining the action
of $Out(F)$ on $Q(F)$ are, in a sense linear. Indeed, for a fixed
$q\in Q(F)$ all the coordinate functions $r_v(q)=\phi(q)_v$, (where
$v\in F$) are fractional-linear functions with the same denominators
(independent of $v$). 
\end{rem}

\section{Applications to the geometry of automorphisms}

In this section we will prove 
Theorem~\ref{A} and Theorem~\ref{B} stated in the introduction.

\begin{defn}
  We say that an outer automorphism $\phi\in Out(F)$ is \emph{strictly
    hyperbolic} if there is $\lambda>1$ such that
\[
\lambda ||w||\le \max\{ ||\phi(w)||, ||\phi^{-1}||\} \text{ for every
  cyclic word } w.
\]
An outer automorphism $\phi\in Out(F)$ is \emph{hyperbolic} if there
is $n>0$ such that $\phi^n$ is strictly hyperbolic.
\end{defn}
For an outer automorphism $\phi\in Out(F)$ set
\[
\lambda_0(\phi):=\inf_{w\ne 1} \, \max\{ \frac{||\phi(w)||}{||w||},
\frac{||\phi^{-1}||}{||w||}\}.
\]

\begin{thm}\label{thm:hyperb}
  The following hold:
\begin{enumerate}
\item An an outer automorphism $\phi\in Out(F)$ is strictly hyperbolic
  if and only if $\lambda_0(\phi)>1$.
\item For any $\phi\in Out(F)$ the number $\lambda_0(\phi)$ is
  rational.
\item There is a double exponential time (in terms of the number $t$
  of Nielsen automorphisms in the expression of $\phi$) algorithm that
  computes $\lambda_0(\phi)$ and decides if an outer automorphism
  $\phi\in Out(F)$ is strictly hyperbolic.
\end{enumerate}
\end{thm}

\begin{proof}
  Let $\phi$ be an outer automorphism given as a product (modulo $Inn
  F$) of $t$ Nielsen automorphisms. For  $u\in F$ with
  $|u|=8^t$ let $d(u,\phi)$ and $d(u,\phi^{-1})$ be the integer constants
  provided by Proposition~\ref{phi1}.

  Then for any cyclic word $w$ by Proposition~\ref{phi1} we have

\[
\frac{||\phi(w)||}{||w||}=\sum_{|u|=8^t} d(u,\phi) f_w(u)
\]
and 
\[
\frac{||\phi^{-1}(w)||}{||w||}=\sum_{|u|=8^t} d(u,\phi^{-1}) f_w(u).
\]

Define the functions $g,h: Q_{8^t} \to\mathbb R$ by formulas
\[
g(p):= \sum_{|u|=8^t} d(u,\phi) p_u \quad \text{ and }
\quad h(p):=\sum_{|u|=8^t} d(u,\phi^{-1}) p_u.
\]

Thus both $h$ and $g$ are affine functions with integer coefficients
on $Q_{8^t}$. Also, put $e:=\max\{g,h\}$. By definition $\phi$ is
strictly hyperbolic if and only if $\inf_{p\in \alpha_{8^t}(C)} e(p)>1$. Since
$e$ is continuous and $\alpha_{8^t}(C)$ is dense in a compact space
$Q_{8^t}$, this condition is equivalent to:
\[
z:=\min_{p\in Q_{8^t}} e(p)>1. \tag{$\heartsuit$}
\]
Moreover, again by density of $\alpha_{8^t}(C)$ in $Q_{8^t}$ if $z>1$
then $z=\lambda_0(\phi)$. Since $h,g$ are affine functions with
integer coefficients on a polyhedron $Q_{8^t}$ given by affine
equations and inequalities with integer coefficients, this implies
that $z$ is a rational number. Thus we have established statements (1)
and (2) of the Theorem.

The number of $v\in F$ with $|u|=8^t$ is $D(8^t)=2k(2k-1)^{8^t}$,
which is double-exponential in $t$. Since a linear programming problem
is solvable in polynomial time, both in the number of variables and
the size of the entries (see, for example,~\cite{Kar}), the value of
$z$ can be computed algorithmically is double exponential time in $t$
and part (3) of the theorem follows.
\end{proof}

Note that Theorem~\ref{thm:hyperb} provides the following new
algorithm for deciding if $\phi\in Out(F)$ is hyperbolic. By
Brinkmann's result~\cite{Br} $\phi$ is hyperbolic if and only if it
does not have nontrivial periodic conjugacy classes. Therefore we will
run the following two procedures in parallel. For $n=1, 2,\dots$ we
start checking, using Theorem~\ref{thm:hyperb}, if $\phi^n$ is
strictly hyperbolic. At the same time for each such $n$ we will start
enumerating elements of $F$ and checking if $\phi^n(f)$ is ever
conjugate to $f$. Eventually we will find $n$ such that either
$\phi^n$ is strictly hyperbolic or $\phi^n$ fixes a nontrivial
conjugacy class.

The only other algorithm for checking hyperbolicity of $\phi$ known to
the author proceeds as follows. Form the group $G=F\rtimes_{\phi}
\mathbb Z$ and in parallel start checking if $G$ is word-hyperbolic
(using, say, Papasoglu's algorithm~\cite{Pa}) and if $\phi$ has a
periodic conjugacy class.  By Brinkmann's theorem one of these
procedures will eventually terminate, and $\phi$ is hyperbolic if and
only if $G$ is word-hyperbolic.

For an outer automorphism $\phi\in Out(F)$ set
\[
\nu_{+}(\phi):=\sup_{w\ne 1} \, \frac{||\phi(w)||}{||w||} \quad\text{
  and }\quad \nu_{-}(\phi):=\inf_{w\ne 1} \,
\frac{||\phi(w)||}{||w||}.
\]

Note that by definition $\frac{1}{\nu_{+}(\phi^{-1})}=\nu_{-}(\phi)$.

The following is a restatement of Theorem~\ref{A}.

\begin{thm}
  Let $\phi\in Out(F)$. Then:
\begin{enumerate}
\item The numbers $\nu_{+}(\phi), \nu_{-}(\phi)$ are rational and
  algorithmically computable in terms of $\phi$.
  
\item For any rational number $\nu_{-}(\phi)<r<\nu_{+}(\phi)$ there
  exists a nontrivial cyclic word $w$ with
\[
\frac{||\phi(w)||}{||w||}=r.
\]

\item We have $\nu_{+}(\phi), \nu_{-}(\phi)\in I(\phi)$.
\end{enumerate}

\end{thm}
\begin{proof}
  
  Let $\phi$ be an outer automorphism given as a product (modulo $Inn
  F$) of $t$ Nielsen automorphisms. Let for $x\in X$ and $u\in F$ with
  $|u|=8^t$ let $d(u,\phi)$ be the constants provided by
  Proposition~\ref{phi1}.
  
  Then as in the proof of Theorem~\ref{thm:hyperb} for any cyclic word
  $w$ we have

\[
\frac{||\phi(w)||}{||w||}= \sum_{|u|=8^t} d(u,\phi)
f_w(u)
\]

As before, define the function $g: Q_{8^t} \to\mathbb R$ as:

\[
g(p):= \sum_{|u|=8^t} d(u,\phi) p_u
\]

Then
\[\nu_{+}(\phi)=\sup_{p\in \alpha_{8^t}(C)} g(p)=\max_{p\in Q_{8^t}} g(p),
\]
where the last equality holds since $\alpha_{8^t}(C)$ is dense in
$Q_{8^t}$.

Similarly,
Then
\[\nu_{-}(\phi)=\inf_{p\in \alpha_{8^t}(C)} g(p)=\min_{p\in Q_{8^t}} g(p).
\]

Recall that $Q_{8^t}$ is a compact convex polyhedron given by
equations and inequalities with rational coefficients and that $g(p)$
is a linear function with integer coefficients. Hence the maximum of
$g(p)$ over $Q_{8^t}$ is a rational number that is algorithmically
computable. The same applies to the minimum of $g$ over $Q_{8^t}$.
This proves assertion (1) of the theorem.

To see that (2) holds, suppose now that $r$ is a rational number satisfying
$\nu_{-}(\phi)<r<\nu_{+}(\phi)$. Let $p_+,p_{-}\in Q_{8^t}$ be such
that $g(p_{+})=\nu_{+}(\phi)$ and $g(p_{-})=\nu_{-}(\phi)$.

Since $QQ_{8^t}^+$ is dense in $Q_{8^t}$, there exist $p_n,q_n\in
QQ_{8^t}^+$ such that $p_n\to p_{+}$ and $q_n\to p_{-}$ as
$n\to\infty$.

Hence $\lim_{n\to\infty} g(p_n)=\nu_{+}(\phi)$ and $\lim_{n\to\infty}
g(q_n)=\nu_{-}(\phi)$. Therefore for some $n$ we have
\[
g(q_n)<r<g(p_n).
\]
Choose a rational number $s\in (0,1)$ so that $r=sg(q_n)+(1-s)g(p_n)$.
Recall that both $p_n,q_n$ are vectors in ${\mathbb R}^{D(8^t)}$ with
all their coordinates being positive rational numbers. Therefore the
same is true for their convex combination $y_n:=sq_n+(1-s)p_n$. Since
$Q_{8^t}$ is convex, $y_n\in Q_{8^t}$ and hence $y_n\in QQ_{8^t}^+$.
By linearity of $g$ we have $g(y_n)=r$. Since by Theorem~\ref{realize}
$QQ_{8^t}^+$ is contained in $\alpha_{8^t}(C)$, there is a cyclic word
$w$ with $\alpha_{8^t}(w)=y_n$. Then
\[
g(y_n)=\frac{||\phi(w)||}{||w||}=r,
\]
as required.

Note that $g: Q_{8^t}\to \mathbb R$ is a linear function on a compact
convex finite-dimensional polyhedron $Q_{8^t}$. Therefore there exists
an extremal point $z$ of $Q_{8^t}$ such that $g(z)=\max_{Q_{8^t}}
g=\nu_{+}(\phi)$. By Lemma~\ref{ext} the improved intial graph
$\Gamma'_z$ is connected and hence by Theorem~\ref{realize} there
exists $w\in C$ such that $z=\alpha_{8^t}(w)$. Then 
\[
\frac{||\phi(w)||}{||w||}=g(z)=\nu_{+}(\phi),
\]
so that $\nu_{+}(\phi)\in I(\phi)$.

The same argument, with replacing $\max$ by $\min$, shows that
$\nu_{-}(\phi)\in I(\phi)$, which completes the proof of the theorem.
\end{proof}

\begin{rem}
One can provide a slightly more uniform argument for Algorithmic computability of $\nu_{\pm}(\phi)$ than that given in the above proof.

In the notation of that proof, a linear function $g$ on the compact convex polyhedron $Q_{8^t}$ attains its maximal and minimal values at some extremal points of $Q_{8^t}$. We have seen in Lemma~\ref{ext} that extremal points of $Q_{8^t}$ have connected improved initial graphs and are thus realizable by some cyclic words in $F$. First, we can find and enumerate all the extremal points of $Q_{8^t}$. For each extremal point $p$, using the Euler circuit algorithm from the proof of Theorem~\ref{realize} we can find a cyclic word $w_p$ realizing that point. We then compute $||\phi(w_p)||/||w_p||$. The maximal of these values 9as $p$ varies over the set of extremal points of $Q_{8^t}$ is $\nu_{+}(\phi)$ and the minimal is $\nu_{-}(\phi)$.

The advantage of this approach is that we do not have to compute an explicit formula for the linear function $g$ at all and the same set of words $w_p$ will work for arbitrary outer automorphisms that are products of at most $t$ Nielsen moves.
\end{rem}

\section{Geodesic Currents}\label{sect:currents}

We provide here a very brief and informal discussion about geodesic
currents on free groups. The details can be found in the dissertation
of Martin~\cite{Ma}. We also intend to elaborate this discussion in
the future. An excellent source on geodesic currents in the context of
Gromov-hyperbolic spaces is provide by a paper of Furman~\cite{Furman}.

Let $F$ be a finitely generated free group. In this section we will
think of $\partial F$ as the topological boundary of $F$ in the sense
of the theory of Gromov-hyperbolic groups. Thus as a set $\partial F$
consists of all equivalence classes of sequences $(g_n)_{n\ge 1}$ of elements $g_n\in F$ that
are \emph{convergent at infinity}. A sequence $(g_n)_n$ is convergent
at infinity if for some (and hence for any) free basis $X$ of $F$ the
distance from $1$ to the geodesic segment from $g_n$ to $g_m$ in the
Cayley graph of $F$ with respect to $X$ tends to infinity as
$m,n\to\infty$. Two sequences $(g_n)_{n\ge 1}$ and $(h_n)_{n\ge 1}$
convergent at infinity are equivalent if the sequence
$g_1,h_1,g_2,h_2,\dots, g_n,h_n,\dots$ is also convergent at
infinity. For each free basis $X$ of $F$ there is a canonical
identification of $\partial F$ with the boundary of $F$ in the sense
of the previous sections, that is with the set of  geodesic
rays starting at $1$ in the Cayley graph of $F$ with respect to $X$.
Every such identification endows $\partial F$ with a topology, that
turns out to be independent of the choice of $X$.

\begin{defn}[Geodesic Currents]
Let $F$ be a finitely generated free group. Denote
\[
\partial^2 F:=\{(s,t): s,t\in \partial F, s\ne t\},
\]
the set of ordered pairs of distinct elements of $\partial F$.

A \emph{geodesic current} on $F$ is an $F$-invariant positive Borel
measure on $\partial^2 F$. 

Two currents $\nu_1, \nu_2$ are \emph{equivalent}, denoted $\nu_1\sim
\nu_2$, if there is $c>0$ such that $\nu_1=c \nu_2$. 

Put 
\[
Curr(F):=\{\nu: \nu \text{ is a nonzero geodesic current on } F \}/\sim
\]
\end{defn}

We will denote the $\sim$-equivalence class of a geodesic current
$\nu$ by $[\nu]$. 
Note that it is often customary to require geodesic currents to be not
just $F$-invariant, but also invariant under the flip-map $\partial^2
F\to \partial^2 F$, $(s,t)\to (t,s)$. However, the two notions are
very close and most arguments work in both contexts.

The set of all geodesic currents $\nu$ on $F$ comes equipped
with the natural weak-* topology. Hence $Curr(F)$ inherits the
quotient topology under the $\sim$-quotient map.

Let $X$ be a free basis of $F$. Then $X$ defines a generating set of
the Borel $\sigma$-algebra of $\partial^2 F$ as follows.
Let $\Gamma(F,X)$ be the Cayley graph of $F$ with respect to $X$.
For any two distinct points $s,t\in F\cup \partial F$ there exists a
unique directed geodesic line $[s,t]_A\subseteq \Gamma(F,X)$ from $s$
to $t$. Thus we think of $[s,t]_X$ as a union of closed edges in
$\Gamma(F,X)$ together with a choice of direction (from $s$ to $t$). 
We will say that
$[s,t]_X\subseteq [s',t']_X$ if the underlying set of $[s,t]_X$ is
contained in that of $[s',t']_X$ and the directions on $[s,t]_X$ and
$[s',t']_X$ agree.

For any $u_1,u_2\in F, u_1\ne u_2$ put 
\[
Cyl_X(u_1,u_2):=\{(s,t)\in \partial^2 F:  [u_1,u_2]_X\subseteq [s,t]_X\} 
\]

It is easy to see that $\mathcal B_X:=\bigcup_{(u,u')\in F^2, u\ne u'} Cyl_X(u,u')$ is
a generating set for the Borel $\sigma$-algebra of $\partial^2 F$.

The following lemma shows that the frequency space $Q_A(F)$ is
essentially the same object as the space of currents $Curr(F)$.

\begin{lem}
Let $A$ be a free basis of $F$. For each $\mu\in Q_X(F)$ define a set
function $\beta_X(\mu): \mathcal B_X\to \mathbb R$ as
$\beta_X(\mu)(Cyl_X(u_1,u_2)):=\mu(Cyl_X(u_1^{-1}u_2))$.

Then for every $\mu\in Q_X(F)$ the function $\beta_X(\mu)$ defines a
geodesic current on $F$ and the map
$\hat\beta_X: Q_X(F)\to Curr(F), \mu\mapsto [\beta_X(\mu)]$ is a homeomorphism.
\end{lem}

\begin{proof}
The proof is a straightforward corollary of the definitions and we
will omit most of the details.

We will indicate why the map  $\hat\beta_X$ is onto. Let $\nu$ be a
geodesic current. After multiplying $\nu$ by an appropriate constant
(which preserves the $\sim$-equivalence class of $\nu$),
we may assume that
\[
\sum_{x\in X\cup X^{-1}} \nu (Cyl_X(1,x))=1.
\]
Since $\nu$ is finitely additive, a simple inductive argument now
shows that for any $m\ge 1$
\[
\sum_{u\in F, |u|_X=m} \nu (Cyl_X(1,u))=1.
\]

For each $u\in F, u\ne 1$ define $\mu(Cyl_X(u)):=\nu (Cyl_X(1,u))$.
We claim that $\mu\in Q_X(F)$.
Indeed, for any $u\in F, u\ne 1$ we have
\[
Cyl_X(1,u)=\sqcup_{|x|=1, |ux|=|u|+1} Cyl_X(1,ux).
\]
Since $\nu$ is a measure, we have
\begin{align*}
&\mu(Cyl_X(u))=\nu(Cyl_X(1,u))=\sum_{|x|=1, |ux|=|u|+1}
\nu(Cyl_X(1,ux) =\\
&=\sum_{|x|=1, |ux|=|u|+1} \mu(Cyl_X(ux).
\end{align*}
This shows that $\mu$ is a Borel probability measure on $\partial F$.

Also, for any $u\in F, u\ne 1$ we have
\[
Cyl_X(1,u)=\sqcup_{|y|=1, |yu|=|u|+1} Cyl_X(y^{-1},u).
\]
Since $\nu$ is $F$-invariant, we have
$\nu(Cyl_X(y^{-1},u))=\nu(Cyl_X(1,yu))$.
Hence we have
\[
\mu(Cyl_X(u))=\sum_{|y|=1, |yu|=|u|+1}\mu(Cyl_X(yu)).
\]
Hence $\mu$ is shift-invariant (with respect to the shift
$T_X$ corresponding to $X$) and hence $\mu\in Q_X(F)$.
By construction we have $\beta_X(\mu)=\nu$, so that $\hat\beta_X$ is onto,
as required.
\end{proof}

We can define an action of $Aut(F)$ on $Curr(F)$ as follows. 
Let $\nu$ be a nonzero geodesic current on $F$ and let $\phi\in
Aut(F)$ be an automorphism. Thus $\phi$ extends to a homeomorphism of
$\partial F$ and so to a homeomorphism of $\partial^2 F$. Define a
positive Borel measure $\phi\nu$ on $\partial^2 F$ by
$(\phi\nu)(B):=\nu(\phi^{-1}B)$ for an arbitrary Borel $B\subseteq
\partial F$. Obviously, $\phi\nu$ is a positive Borel measure on
$\partial^2 F$. It can also be verified directly, via a
Gromov-hyperbolic argument, that $\phi\nu$ is $F$-invariant, and hence
is a geodesic current. 

The measure $\phi\nu$ can be understood ``coordinatewise'' as
follows. Let $X$ be a free basis of $F$, so that $X'=\phi(X)$ is also
a free basis. Then $\phi\nu$ has the same $\mathcal B_{X'}$
``coordinates'' as the $\mathcal B_X$-coordinates of $\nu$. That is,
for any $u_1,u_2\in F, u_1\ne u_2$ we have
$(\phi\nu)(Cyl_{X'}(\phi(u_1),\phi(u_2))=\nu(Cyl_X(u_1,u_2))$.
This defines a continuous action of $Aut(F)$ on $Curr(F)$. Moreover,
if $\alpha$ is an inner automorphism of $F$ corresponding to
conjugation by $h$, $g\mapsto h gh^{-1}$, then $\alpha(t)=ht$ for any
$t\in \partial F$ and hence $\alpha B=hB$ for any $B\subseteq \partial^2
F$. This easily implies that $\phi\nu$ depends only on the outer
automorphism class of $\phi$ in $Out(F)=Aut(F)/Inn(F)$ and that the
action of $Aut(F)$ on $Curr(F)$ factors through to an action of
$Out(F)$ on $Curr(F)$.

If we fix a free basis $X$ of $F$, then one can show that this action
of $Out(F)$ on $Curr(F)$ is exactly the same (via the
$\beta_X$-identification of $Curr(F)$ and $Q_X(F)$ as the action of
$Out(F)$ on $Q_X(F)$ constructed earlier. One way to see this is to check
that the actions coincide on the dense set of measures coming from the
conjugacy classes in $F$, that is $\alpha_X(C)\subseteq Q_X(F)$, as
well as the corresponding dense set of geodesic currents
$\hat\beta_X\alpha_X(C)\subseteq Curr(F)$.

There is also a natural topological embedding $j:{\mathcal X}(F)\to Curr(F)$ of
the Culler-Vogtmann~\cite{CV}outer space ${\mathcal X}(F)$ into the space of currents $Curr(F)$
(and hence into $Q_X(F)$ as well). A point $\omega$ of an outer space can be
represented by a minimal  free discrete isometric action of $F$ on an
$\mathbb R$-tree $Y$ (or the ``marked length spectrum'' on $F$ defined
by such an action). This tree $Y$ can be thought of as the universal
cover of a finite metric graph $K$. The tree $Y$ is a
Gromov-hyperbolic metric space and as such $\partial Y$ comes equipped
with the so-called Patterson-Sullivan measure. This
measure, via an explicit formula involving the Busemann function, defines an $F$-invariant measure on
$\partial^2 Y$ and hence, via the orbit map $F\to Y$, on $\partial F$. The $\sim$-equivalence class
of that measure is precisely $j(\omega)$.  We
refer the reader to an article by Alex Furman~\cite{Furman} for a
detailed discussion about Patterson-Sullivan measures in the context
of Gromov-hyperbolic spaces and the way in which they determine
geodesic currents.

 Reiner Martin~\cite{Ma} produced a different family of $Out(F_n)$-equivariant embeddings of the outer space into the space of currents, which, since $Curr(F)$ is compact, provide ways of equivariantly compactifying the outer space. Martin~\cite{Ma} proved, that these
compactifications are different from the standard one (given by weak
limits of length functions).

 Yet another way of embedding the outer space into $Curr(F)$ comes from the so-called "visibility currents" defined by weighted non-backtracking random walks on the tree $Y=\widetilde K$ associated with a point of the outer space in the above notations. 

As demonstrated by the work of Lyons~\cite{Ly} for simplicial graphs, one should expect the Patterson-Sullivan embedding, the "visibility current embedding" and Reiner Martin's embeddings of the outer space in the space of currents to be rather different. 
The precise relationship between these  embeddings and the properties of the compactifications of the outer space  that they provide remain to be explored.


\begin{thebibliography}{ABC}
  
\bibitem{BFH} M.~Bestvina, M.~Feighn and M.~Handel, \emph{Laminations
    and irreducible automorphisms of free groups}, Geometric and
  Funct.  Anal. \textbf{7} (1997), no. 2, 215--244

\bibitem{BF92}
M.~Bestvina and M.~Feighn, \emph{The combination theorem for
negatively curved groups}, J. Diff. Geom. \textbf{35} (1992),
85--101

\bibitem{BF96}
M.~Bestvina and M.~Feighn, \emph{Addendum and correction to: "A
combination theorem for negatively curved groups",}  J. Diff.
Geom. \textbf{43} (1996), no. 4, 783--788

\bibitem{BF00}
M.~Bestvina and M.~Feighn, \emph{The topology at infinity of $Out(F_n)$.}  Invent. Math.  140  (2000),  no. 3, 651--692

\bibitem{BH92}
M.~Bestvina, and M.~Handel, \emph{Train tracks and automorphisms of free groups.}  Ann. of Math. (2)  \textbf{135}  (1992),  no. 1, 1--51


\bibitem{Bo86}
F.~Bonahon, \emph{Bouts des vari\'et\'es hyperboliques de
dimension $3$.}  Ann. of Math. (2) \textbf{124} (1986), no. 1,
71--158

\bibitem{Bo88}
F.~Bonahon, \emph{The geometry of Teichmüller space via geodesic
currents.} Invent. Math. \textbf{92} (1988), no. 1, 139--162

\bibitem{Bo91}
F.~Bonahon, \emph{Geodesic currents on negatively curved groups.}
Arboreal group theory (Berkeley, CA, 1988), 143--168, Math. Sci.
Res. Inst. Publ., 19, Springer, New York, 1991

\bibitem{BV}
M.~Bridson, and K.~Vogtmann, \emph{The symmetries of outer space.}  Duke Math. J.  \textbf{106}  (2001),  no. 2, 391--409.


\bibitem{Br}
P.~Brinkmann, \emph{Hyperbolic automorphisms of free groups,}
Geometric and Functional Analysis,  \textbf{10}  (2000),  no. 5, 1071--1089

\bibitem{CV}
M.~Culler, K.~Vogtmann, \emph{Moduli of graphs and automorphisms of free groups.}  Invent. Math.  \textbf{84}  (1986),  no. 1, 91--119.


\bibitem{Furman}
A.~Furman, \emph{Coarse-geometric perspective on negatively curved
  manifolds and groups,} 
in ``Rigidity in Danamics and Geometry (editors M. Burger and
A. Iozzi)'', Springer 2001, 149--166


\bibitem{Ga}
F.~Gautero, \emph{Hyperbolicity of mapping class groups and
spaces}, Enseign. Math. (2)  \textbf{49}  (2003),  no. 3-4, 263--305


\bibitem{KKS} V.~Kaimanovich, I.~Kapovich and P.~Schupp,
\emph{Generic stretching factors of group homomorphisms,} in preparation


\bibitem{KSS} I.~Kapovich, P.~Schupp, and V.~Shpilrain, \emph{Generic
    properties of Whitehead's algorithm and isomorphism rigidity of
    random one-relator groups}, Pacific J. Math., to appear


\bibitem{Kar}
H.~Karloff, \emph{Linear programming.} Progress in Theoretical Computer Science. Birkhäuser Boston, Inc., Boston, MA, 1991

\bibitem{KH} A.~Katok, and B.~Hasselblatt, \emph{Introduction to the modern theory of dynamical systems.} With a supplementary chapter by Katok and Leonardo Mendoza. Encyclopedia of Mathematics and its Applications, 54. Cambridge University Press, Cambridge, 1995

\bibitem{LL}
G.~Levitt, and M.~Lustig, \emph{Irreducible automorphisms of $F_n$ have north-south dynamics on compactified outer space.}  J. Inst. Math. Jussieu  \textbf{2}  (2003),  no. 1, 59--72.

\bibitem{LS} R.~Lyndon and P.~Schupp, \emph{Combinatorial Group
    Theory,} Springer-Verlag, 1977. Reprinted in the ``Classics in
  mathematics'' series, 2000.


\bibitem{Ly}
R.~Lyons, 
\emph{Equivalence of boundary measures on covering trees of finite graphs,}
Ergodic Theory Dynam. Systems \textbf{14} (1994), no. 3, 575--597

\bibitem{Ma} R.~Martin, \emph{Non-Uniquely Ergodic Foliations of Thin
    Type, Measured Currents and Automorphisms of Free Groups}, PhD
  Thesis, 1995


\bibitem{Ox}
J.~Oxtoby, \emph{On two theorems of Parthasarathy and Kakutani concerning the shift transformation.} 1963  Ergodic Theory (Proc. Internat. Sympos., Tulane Univ., New Orleans, La., 1961)  pp. 203--215 Academic Press, New York

\bibitem{Pa} P.~Papasoglu, \emph{An algorithm detecting hyperbolicity.}  Geometric and computational perspectives on infinite groups,  193--200, DIMACS Ser. Discrete Math. Theoret. Comput. Sci., 25, Amer. Math. Soc., Providence, RI, 1996

\bibitem{Par}
K.~R.~Parthasarathy, \emph{On the category of ergodic measures.}  Illinois J. Math.  \textbf{5}  (1961), 648--656

\bibitem{West} D.~West, \emph{Introduction to Graph Theory,} (2-nd
  edition) Prentice Hall, 2001

\end{thebibliography}
\end{document}